# Some Consequences of the Phase Space Choice

**Stefan Balint[a] and Agneta M. Balint[b]**

[a]*West University of Timisoara, Department of Computer Science, Bulv. V.Parvan 4, 300223 Timisoara, Romania*
[b]*West University of Timisoara, Faculty of Physics, Bulv. V.Parvan 4, 300223 Timisoara, Romania*

*Corresponding author*: Stefan Balint, *West University of Timisoara, Faculty of Mathematics and Computer Science, Bulv. V.Parvan 4, 300223 Timisoara, Romania*
E-mail: stefan.balint@e-uvt.ro

**Abstract.** In the literature, for semidynamical systems in infinite dimensional phase spaces, different topological structures are used (Hilbert, Banach, Sobolev, locally convex, Hausdorf topology etc.). That is because there are neither set rules nor understanding of the "right way" to choose the phase space and its topology starting from a system of partial differential equations. The goal of this paper is to reveal the influence of the choice of the phase space and its topology as concern the results obtained for the semidynamical system defined by the same system of partial differential equations. In the paper the linear 3D Euler equations are considered which are obtained by linearizing the non linear 3D Euler equations at a constant solution. The well posedness of the instantaneous perturbation propagation problem, that of the permanent source produced time harmonic perturbation propagation problem as well the stability of the null solution are analyzed in three different phase spaces. The idea is to derive explicit solutions for the linear 3D Euler equations, to build up different phase spaces by using explicit solutions and analyze the well-posedness as well the stability of the null solution. The obtained results present significant differences and some of them are surprising. For instance, in one phase space the stability of the null solution coexists with solutions having strictly positive exponential growth rate, but in other one the propagation problem is ill posed. This aspect is extremely important from the point of view of a real phenomenon modeled by the linear 3D Euler equations.

## 1. Introduction

In the 3D gas flow model, the nonlinear Euler equations governing the flow of an inviscid, compressible, non heat conducting, isentropic, perfect gas are [10]:

$$\frac{\partial v_x}{\partial t} + v_x \cdot \frac{\partial v_x}{\partial x} + v_y \cdot \frac{\partial v_x}{\partial y} + v_z \cdot \frac{\partial v_x}{\partial z} + \frac{1}{\rho}\frac{\partial p}{\partial x} = 0$$

$$\frac{\partial v_y}{\partial t} + v_x \cdot \frac{\partial v_y}{\partial x} + v_y \cdot \frac{\partial v_y}{\partial y} + v_z \cdot \frac{\partial v_y}{\partial z} + \frac{1}{\rho}\frac{\partial p}{\partial y} = 0 \qquad (1)$$

$$\frac{\partial v_z}{\partial t} + v_x \cdot \frac{\partial v_z}{\partial x} + v_y \cdot \frac{\partial v_z}{\partial y} + v_z \cdot \frac{\partial v_z}{\partial z} + \frac{1}{\rho} \frac{\partial p}{\partial z} = 0$$

$$\frac{\partial \rho}{\partial t} + v_x \cdot \frac{\partial \rho}{\partial x} + v_y \cdot \frac{\partial \rho}{\partial y} + v_z \cdot \frac{\partial \rho}{\partial z} + \rho \cdot \left( \frac{\partial v_x}{\partial x} + \frac{\partial v_y}{\partial y} + \frac{\partial v_z}{\partial z} \right) = 0.$$

Here: $t$ - time, $v_x$, $v_y$, $v_z$ - velocity components along the $Ox, Oy, Oz$ axis respectively; $p$ - pressure, $\rho$ - density.

The pressure $p$, the density $\rho$ and the absolute temperature $T$ satisfy the state equation of the perfect gas

$$p = \rho \cdot R \cdot T \tag{2}$$

where the gas constant $R = c_p - c_v$; $c_p$ and $c_v$ being the specific heat capacities at constant pressure and constant volume, respectively.

If $v_x = U_0 = const > 0$; $v_y = 0$; $v_z = 0$; $\rho = \rho_0 = const > 0$; $p = p_0 = const > 0$, then according to (2): $p_0 = \rho_0 \cdot R \cdot T_0$ and the associated isentropic sound speed $c_0$ verifies $c_0^2 = \frac{p_0}{\rho_0} = R \cdot T_0$.

Linearizing (1) at $v_x = U_0$; $v_y = v_z = 0$; $p = p_0$; $\rho = \rho_0$ and using that the perturbations $p'$, $\rho'$ satisfy:

$$\left( \frac{\partial}{\partial t} + U_0 \cdot \frac{\partial}{\partial x} \right)\left( p' - c_0^2 \cdot \rho' \right) = 0 \tag{3}$$

the following system of linear Euler equations is obtained:

$$\frac{\partial v'_x}{\partial t} + U_0 \cdot \frac{\partial v'_x}{\partial x} + \frac{1}{\rho_0} \cdot \frac{\partial p'}{\partial x} = 0$$

$$\frac{\partial v'_y}{\partial t} + U_0 \cdot \frac{\partial v'_y}{\partial x} + \frac{1}{\rho_0} \cdot \frac{\partial p'}{\partial y} = 0 \tag{4}$$

$$\frac{\partial v'_z}{\partial t} + U_0 \cdot \frac{\partial v'_z}{\partial x} + \frac{1}{\rho_0} \cdot \frac{\partial p'}{\partial z} = 0$$

$$\frac{\partial p'}{\partial t} + U_0 \cdot \frac{\partial p'}{\partial x} + \rho_0 \cdot c_0^2 \left( \frac{\partial v'_x}{\partial x} + \frac{\partial v'_y}{\partial y} + \frac{\partial v'_z}{\partial z} \right) = 0.$$

**DEFINITION 1.** Following [9], for a given set $\mathcal{J} = \{I\}$ of initial values $I = (F, G, H, P)$ one calls the initial value problem (4), (5)

$$v'_x(x, y, z, 0) = F; \quad v'_y(x, y, z, 0) = G; \quad v'_z(x, y, z, 0) = H; \quad p'(x, y, z, 0) = P \tag{5}$$

well-posed on the interval of time $[0, t_0]$ if:

a) for any $I \in \mathcal{J}$ there is a unique solution to the initial value problem (4), (5), defined for any $t \in [0, t_0]$;

b) the solution varies continuously with the initial data $I$ (i.e. for any $I \in \mathcal{J}$ and any sequence $(I_n)_n$, $I_n \in \mathcal{J}$ if $I_n$ converges to $I$ in $\mathcal{J}$, then the sequence of solutions corresponding to $I_n$, converges in $\mathcal{J}$ to the solution corresponding to $I$ at any $t \in [0, t_0]$).

In this paper the initial value problem (4), (5) is called also ***instantaneous perturbation propagation problem.***
Beside the instantaneous perturbation propagation problem (4), (5) we will consider also the so called permanent source produced time harmonic perturbation propagation problem.

**DEFINITION 2**. The solution of the system of non-homogeneous linear Euler equations:

$$\frac{\partial v'_x}{\partial t} + U_0 \cdot \frac{\partial v'_x}{\partial x} + \frac{1}{\rho_0} \cdot \frac{\partial p'}{\partial x} = h(t) \cdot F \cdot \sin \omega_f t$$

$$\frac{\partial v'_y}{\partial t} + U_0 \cdot \frac{\partial v'_y}{\partial x} + \frac{1}{\rho_0} \cdot \frac{\partial p'}{\partial y} = h(t) \cdot G \cdot \sin \omega_f t \qquad (6)$$

$$\frac{\partial v'_z}{\partial t} + U_0 \cdot \frac{\partial v'_z}{\partial x} + \frac{1}{\rho_0} \cdot \frac{\partial p'}{\partial z} = h(t) \cdot H \cdot \sin \omega_f t$$

$$\frac{\partial p'}{\partial t} + U_0 \cdot \frac{\partial p'}{\partial x} + \rho_0 \cdot c_0^2 \left( \frac{\partial v'_x}{\partial x} + \frac{\partial v'_y}{\partial y} + \frac{\partial v'_z}{\partial z} \right) = h(t) \cdot P \cdot \sin \omega_f t$$

which is equal to zero for $t \leq 0$ is called ***solution of the permanent source produced time harmonic perturbation propagation problem.***
In equations (6) $A = (F, G, H, P)$ is the amplitude of the permanent source produced time harmonic perturbation; $\omega_f$ is the angular frequency and $h(t)$ is the Heaviside function.
$F, G, H, P$ are real valued functions depending on the spatial variables $x, y, z$, like in equations (5).

**DEFINITION 3**. Following [9], for a given class $\mathcal{A} = \{A\}$ of amplitudes $A = (F, G, H, P)$ one calls the permanent source produced time harmonic perturbation propagation problem (6) ***well-posed*** on the interval of time $[0, t_0]$ if:
a) for any $A \in \mathcal{A}$ there is a unique solution to the problem (6) defined for any $t \in [0, t_0]$
b) the solution varies continuously with amplitude $A$ (i.e. for any $A \in \mathcal{A}$ and any sequence $(A_n)_n$, $A_n \in \mathcal{A}$ if $A_n$ converges to $A$ in $\mathcal{A}$, then the sequence of solutions corresponding to $A_n$ converges in $\mathcal{A}$ to the solution corresponding to $A$ at any $t \in [0, t_0]$.

**DEFINITION 4**. The ***null solution*** of (4) ***is stable*** (in sense of Lyapunov) with respect to the class of the initial value perturbations $\mathcal{J}$ (to the set of permanent source produced time harmonic perturbations whose amplitude is in the class $\mathcal{A}$, respectively) if:
a) the instantaneous perturbation propagation problem (4), (5) (the permanent source produced time harmonic perturbation propagation problem (6)) is well-posed on the interval of time $[0, t_0]$ for any $t_0 > 0$;

b) the solution of the instantaneous perturbation propagation problem (4), (5) (the solution of the permanent source produced time harmonic perturbation propagation problem (6)) is close to the null solution of (4) all time $t \geq 0$ provided it is sufficiently close at $t = 0$ (it is close to the null solution of (4) all time provided the amplitude $A$ is sufficiently close to zero).

The precise meaning of the concepts "perturbation", "small perturbation", "solution to the perturbation propagation problem", "small solution" has to be designed by other definitions and there is some freedom here. Using this freedom, in the following we present different phase spaces in order to reveal that the mathematical results and tools are specific.

## 2. Explicit Solutions

The next statement concerns the existence of a class of explicit point-wise solutions of the instantaneous perturbation propagation problem (4), (5).

**PROPOSITION 5.** *For any system* $f_i$, $i = \overline{1,4}$ *of continuously differentiable functions* $f_i : R^1 \to R^1$ *and any system* $k_i, l_i, m_i$, $i = \overline{1,4}$ *of real numbers with* $k_3 \cdot k_4 \neq 0$ *the functions defined by:*

$$v_x'(x,y,z,t) = k_1 \cdot f_1\left[k_1 x + l_1 y + m_1 z - \left(k_1 U_0 - c_0 \sqrt{k_1^2 + l_1^2 + m_1^2}\right) \cdot t\right] +$$

$$k_2 \cdot f_2\left[k_2 x + l_2 y + m_2 z - \left(k_2 U_0 + c_0 \sqrt{k_2^2 + l_2^2 + m_2^2}\right) \cdot t\right] +$$

$$\frac{l_3}{k_3} \cdot f_3[k_3 x + l_3 y + m_3 z - k_3 U_0 t] + \frac{m_4}{k_4} \cdot f_4[k_4 x + l_4 y + m_4 z - k_4 U_0 t]$$

$$v_y'(x,y,z,t) = l_1 \cdot f_1\left[k_1 x + l_1 y + m_1 z - \left(k_1 U_0 - c_0 \sqrt{k_1^2 + l_1^2 + m_1^2}\right) \cdot t\right] +$$

$$l_2 \cdot f_2\left[k_2 x + l_2 y + m_2 z - \left(k_2 U_0 + c_0 \sqrt{k_2^2 + l_2^2 + m_2^2}\right) \cdot t\right] - \qquad (7)$$

$$f_3[k_3 x + l_3 y + m_3 z - k_3 U_0 t]$$

$$v_z'(x,y,z,t) = m_1 \cdot f_1\left[k_1 x + l_1 y + m_1 z - \left(k_1 U_0 - c_0 \sqrt{k_1^2 + l_1^2 + m_1^2}\right) \cdot t\right] +$$

$$m_2 \cdot f_2\left[k_2 x + l_2 y + m_2 z - \left(k_2 U_0 + c_0 \sqrt{k_2^2 + l_2^2 + m_2^2}\right) \cdot t\right] -$$

$$f_4[k_4 x + l_4 y + m_4 z - k_4 U_0 t]$$

$$p'(x,y,z,t) = -c_0 \rho_0 \sqrt{k_1^2 + l_1^2 + m_1^2} \cdot f_1\left[k_1 x + l_1 y + m_1 z - \left(k_1 U_0 - c_0 \sqrt{k_1^2 + l_1^2 + m_1^2}\right) \cdot t\right] +$$

$$c_0 \rho_0 \sqrt{k_2^2 + l_2^2 + m_2^2} \cdot f_2\left[k_2 x + l_2 y + m_2 z - \left(k_2 U_0 + c_0 \sqrt{k_2^2 + l_2^2 + m_2^2}\right) \cdot t\right]$$

*verify point-wise the equations (4) and the initial condition:*
$$v_x'(x,y,z,0) = F(x,y,z) = k_1 \cdot f_1(k_1 x + l_1 y + m_1 z) + k_2 \cdot f_2(k_2 x + l_2 y + m_2 z) +$$

$$\frac{l_3}{k_3} \cdot f_3(k_3 x + l_3 y + m_3 z) + \frac{m_4}{k_4} \cdot f_4(k_4 x + l_4 y + m_4 z)$$

$$v_y'(x,y,z,0) = G(x,y,z) = l_1 \cdot f_1(k_1 x + l_1 y + m_1 z) + l_2 \cdot f_2(k_2 x + l_2 y + m_2 z) - f_3(k_3 x + l_3 y + m_3 z) \quad (8)$$

$$v_z'(x,y,z,0) = H(x,y,z) = m_1 \cdot f_1(k_1 x + l_1 y + m_1 z) + m_2 \cdot f_2(k_2 x + l_2 y + m_2 z) - f_4(k_4 x + l_4 y + m_4 z)$$

$$p'(x,y,z,0) = P(x,y,z) = -c_0 \rho_0 \sqrt{k_1^2 + l_1^2 + m_1^2} \cdot f_1(k_1 x + l_1 y + m_1 z) +$$

$$c_0 \rho_0 \sqrt{k_2^2 + l_2^2 + m_2^2} \cdot f_2(k_2 x + l_2 y + m_2 z)$$

*Proof.* By calculus. □

The following statement reveals that in the class of explicit solutions given by (7) there exist solutions for which velocity potential doesn't exist.

**PROPOSITION 6.** *If* $f_1 = f_2 = f_3 = 0$; $f_4(\xi) = \sin \xi$; $k_4 = l_4 = m_4 = 1$, *then for the explicit point-wise solution (7) velocity potential doesn't exist.*

*Proof.* In the considered case the explicit solution is:
$$v_x'(x,y,z,t) = \sin(x + y + z - U_0 t); \quad v_y'(x,y,z,t) = 0; \quad v_z'(x,y,z,t) = -\sin(x + y + z - U_0 t);$$
$$p'(x,y,z,t) = 0.$$

It follows that: $\dfrac{\partial v_x'}{\partial y} = \cos(x + y + z - U_0 t)$ and $\dfrac{\partial v_y'}{\partial x} = 0$. So, $\dfrac{\partial v_x'}{\partial y} \neq \dfrac{\partial v_y'}{\partial x}$ and there is no function $\Phi(x,y,z,t)$ having the property: $v_x' = \dfrac{\partial \Phi}{\partial x}$; $v_y' = \dfrac{\partial \Phi}{\partial y}$; $v_z' = \dfrac{\partial \Phi}{\partial z}$. □

*Remark.* The class of explicit functions given by (7) can not be derived from a velocity potential $\Phi(x,y,z,t)$.

The next statement concerns the existence of a class of explicit point-wise solutions for the permanent source produced time harmonic perturbation propagation problem (6).

**PROPOSITION 7**. *For any system* $f_i$, $i = \overline{1,4}$ *of continuously differentiable functions* $f_i : R' \to R'$ *and any system* $k_i, l_i, m_i$, $i = \overline{1,4}$ *of real numbers with* $k_3 \cdot k_4 \neq 0$ *the functions defined by:*

$$v_x'(x,y,z,t) = h(t) \cdot k_1 \cdot \int_0^t f_1\left[k_1 x + l_1 y + m_1 z - \left(k_1 U_0 - c_0 \sqrt{k_1^2 + l_1^2 + m_1^2}\right) \cdot (t - \tau)\right] \sin \omega_f \tau \, d\tau +$$

$$h(t) \cdot k_2 \cdot \int_0^t f_2\left[k_2 x + l_2 y + m_2 z - \left(k_2 U_0 + c_0 \sqrt{k_2^2 + l_2^2 + m_2^2}\right) \cdot (t - \tau)\right] \sin \omega_f \tau \, d\tau +$$

$$h(t) \cdot \frac{l_3}{k_3} \cdot \int_0^t f_3[k_3 x + l_3 y + m_3 z - k_3 U_0 (t - \tau)] \sin \omega_f \tau \, d\tau +$$

$$h(t) \cdot \frac{m_4}{k_4} \cdot \int_0^t f_4[k_4 x + l_4 y + m_4 z - k_4 U_0 (t - \tau)] \sin \omega_f \tau \, d\tau$$

$$v_y'(x,y,z,t) = h(t) \cdot l_1 \cdot \int_0^t f_1\left[k_1 x + l_1 y + m_1 z - \left(k_1 U_0 - c_0\sqrt{k_1^2 + l_1^2 + m_1^2}\right) \cdot (t-\tau)\right] \sin \omega_f \tau \, d\tau +$$

$$h(t) \cdot l_2 \cdot \int_0^t f_2\left[k_2 x + l_2 y + m_2 z - \left(k_2 U_0 + c_0\sqrt{k_2^2 + l_2^2 + m_2^2}\right) \cdot (t-\tau)\right] \sin \omega_f \tau \, d\tau -$$

$$h(t) \cdot \int_0^t f_3[k_3 x + l_3 y + m_3 z - k_3 U_0 (t-\tau)] \sin \omega_f \tau \, d\tau$$

(9)

$$v_z'(x,y,z,t) = h(t) \cdot m_1 \cdot \int_0^t f_1\left[k_1 x + l_1 y + m_1 z - \left(k_1 U_0 - c_0\sqrt{k_1^2 + l_1^2 + m_1^2}\right) \cdot (t-\tau)\right] \sin \omega_f \tau \, d\tau +$$

$$h(t) \cdot m_2 \cdot \int_0^t f_2\left[k_2 x + l_2 y + m_2 z - \left(k_2 U_0 + c_0\sqrt{k_2^2 + l_2^2 + m_2^2}\right) \cdot (t-\tau)\right] \sin \omega_f \tau \, d\tau -$$

$$h(t) \cdot \int_0^t f_4[k_4 x + l_4 y + m_4 z - k_4 U_0 (t-\tau)] \sin \omega_f \tau \, d\tau$$

$$p'(x,y,z,t) =$$

$$-h(t) c_0 \rho_0 \sqrt{k_1^2 + l_1^2 + m_1^2} \cdot \int_0^t f_1\left[k_1 x + l_1 y + m_1 z - \left(k_1 U_0 - \sqrt{k_1^2 + l_1^2 + m_1^2}\right)(t-\tau)\right] \sin \omega_f \tau \, d\tau +$$

$$h(t) c_0 \rho_0 \sqrt{k_2^2 + l_2^2 + m_2^2} \int_0^t f_2\left[k_2 x + l_2 y + m_2 z - \left(k_2 U_0 + c_0\sqrt{k_2^2 + l_2^2 + m_2^2}\right)(t-\tau)\right] \sin \omega_f \tau \, d\tau$$

*verify point-wise the equation (6) for:*

$$F = k_1 \cdot f_1(k_1 x + l_1 y + m_1 z) + k_2 \cdot f_2(k_2 x + l_2 y + m_2 z) + \frac{l_3}{k_3} \cdot f_3(k_3 x + l_3 y + m_3 z) +$$

$$\frac{m_4}{k_4} \cdot f_4(k_4 x + l_4 y + m_4 z)$$

$$G = l_1 \cdot f_1(k_1 x + l_1 y + m_1 z) + l_2 \cdot f_2(k_2 x + l_2 y + m_2 z) - f_3(k_3 x + l_3 y + m_3 z)$$

(10)

$$H = m_1 \cdot f_1(k_1 x + l_1 y + m_1 z) + m_2 \cdot f_2(k_2 x + l_2 y + m_2 z) - f_4(k_4 x + l_4 y + m_4 z)$$

$$P = -c_0 \rho_0 \sqrt{k_1^2 + l_1^2 + m_1^2} \cdot f_1(k_1 x + l_1 y + m_1 z) + c_0 \rho_0 \sqrt{k_2^2 + l_2^2 + m_2^2} \cdot f_2(k_2 x + l_2 y + m_2 z)$$

*Proof.* By calculus. □

The next statement reveals that in the class of solutions given by (10) there exist solutions for which velocity potential doesn't exist.

**PROPOSITION 8**. *If $f_1 = f_2 = f_3 = 0$; $f_4(\xi) = \sin \xi$; $k_4 = l_4 = m_4 = 1$, then for the point-wise solution given by (9) velocity potential doesn't exist.*

*Proof.* In the considered case the solution is:

$$v_x'(x,y,z,t) = h(t) \cdot \int_0^t \sin[x + y + z - U_0(t-\tau)] \cdot \sin \omega_f \tau \, d\tau \, ; \quad v_y'(x,y,z,t) = 0 \, ;$$

$$v_z'(x,y,z,t) = -h(t) \cdot \int_0^t \sin[x + y + z - U_0(t-\tau)] \cdot \sin \omega_f \tau \, d\tau \, ; \quad p'(x,y,z,t) = 0.$$

It follows that $\dfrac{\partial v_x'}{\partial y} = h(t) \cdot \int_0^t \cos[x + y + z - U_0(t - \tau)] \cdot \sin \omega_f \tau \, d\tau$; $\dfrac{\partial v_y'}{\partial x} = 0$ and so $\dfrac{\partial v_x'}{\partial y} \neq \dfrac{\partial v_y'}{\partial x}$.

Therefore, there is no function $\Phi(x, y, z, t)$ having the property $v_x' = \dfrac{\partial \Phi}{\partial x}$ and $v_y' = \dfrac{\partial \Phi}{\partial y}$. □

***Remark.*** The class of solutions given by (9) cannot be derived from a velocity potential $\Phi(x, y, z, t)$.

Finally, it has to be noted that the explicit solutions (7) and (9) were found using the method of characteristics described in [4].

## 3. First Phase Space

Let $X_1$ be the topological function space [6], [11] of the set of systems $I = (F, G.H, P)$ (or $A = (F, G, H, P)$) of functions $F, G, H, P : R^3 \to R^1$ which are continuously differentiable, endowed with the usual algebraic operations and topology generated by the uniform convergence on $R^3$ [3]. A neighborhood of the origin $O$ in $X_1$ is a set $V_0$ of systems $I$ from $X_1$ having the property that there exists $\varepsilon \in R^1$, $\varepsilon > 0$, such that if for $I = (F, G.H, P) \in X_1$ inequalities $|F(x, y, z)| < \varepsilon$; $|G(x, y, z)| < \varepsilon$; $|H(x, y, z)| < \varepsilon$; $|P(x, y, z)| < \varepsilon$ hold for any $(x, y, z) \in R^3$, then $I \in V_0$.

The set $V_0^\varepsilon$ defined by:
$$V_0^\varepsilon = \{(F, G, H, P) \in X_1 : |F(x, y, z)| < \varepsilon \text{ and so on, for any } (x, y, z) \in R^3\} \quad (11)$$
is a neighborhood of the origin $O$ in $X_1$.

Let $Y_1$ be the topological function space of the set of systems $f = (f_1, f_2, f_3, f_4)$ of functions $f_i : R^1 \to R^1$, $i = \overline{1,4}$ which are continuously differentiable endowed with the usual algebraic operations and topology generated by the uniform convergence on $R^1$. A neighborhood of the origin $O$ in $Y_1$ is a set $W_0$ of systems $f$ from $Y_1$ having the property that there exists $\varepsilon \in R^1$, $\varepsilon > 0$, such that if for $f = (f_1, f_2, f_3, f_4) \in Y_1$ inequalities $|f_i(\xi)| < \varepsilon$, $i = \overline{1,4}$ hold for any $\xi \in R^1$, then $f \in W_0$.

The set $W_0^\varepsilon$ defined by:
$$W_0^\varepsilon = \{(f_1, f_2, f_3, f_4) \in Y_1 : |f_i(\xi)| < \varepsilon \text{ for any } \xi \in R^1 \text{ and } i = \overline{1,4}\} \quad (12)$$
is a neighborhood of the origin $O$ in $Y_1$.

The next statement concerns the neighborhoods $W_0$, $V_0$ of the origin $O$ in $Y_1$ and $X_1$, respectively.

**PROPOSITION 9**. *For $f = (f_1, f_2, f_3, f_4) \in Y_1$ and real constants $k_i, l_i, m_i$ $i = \overline{1,4}$ with $k_3 \cdot k_4 \neq 0$ and $F, G, H, P$ given by:*

$$F(x,y,z) = k_1 \cdot f_1(k_1 x + l_1 y + m_1 z) + k_2 \cdot f_2(k_2 x + l_2 y + m_2 z) + \frac{l_3}{k_3} \cdot f_3(k_3 x + l_3 y + m_3 z) +$$

$$\frac{m_4}{k_4} \cdot f_4(k_4 x + l_4 y + m_4 z)$$

$$G(x,y,z) = l_1 \cdot f_1(k_1 x + l_1 y + m_1 z) + l_2 \cdot f_2(k_2 x + l_2 y + m_2 z) - f_3(k_3 x + l_3 y + m_3 z) \quad (13)$$

$$H(x,y,z) = m_1 \cdot f_1(k_1 x + l_1 y + m_1 z) + m_2 \cdot f_2(k_2 x + l_2 y + m_2 z) - f_4(k_4 x + l_4 y + m_4 z)$$

$$P(x,y,z) = -c_0 \rho_0 \sqrt{k_1^2 + l_1^2 + m_1^2} \cdot f_1(k_1 x + l_1 y + m_1 z) + c_0 \rho_0 \sqrt{k_2^2 + l_2^2 + m_2^2} \cdot f_2(k_2 x + l_2 y + m_2 z)$$

*the following statements hold:*

i) $I = (F, G, H, P) \in X_1$

ii) if $f \in W_0^\varepsilon$, then $I \in V_0^{\varepsilon \cdot m}$, where $m$ is given by

$$m = \max\{|k_1| + |k_2| + |l_3/k_3| + |m_4/k_4|; |l_1| + |l_2| + 1; |m_1| + |m_2| + 1;$$

$$c_0 \rho_0 \left( \sqrt{k_1^2 + l_1^2 + m_1^2} + \sqrt{k_2^2 + l_2^2 + m_2^2} \right)\} \quad (14)$$

*Proof.* The functions $F, G, H, P$ given by (13) can be evaluated as follows:

$$|F(x,y,z)| \le \left(|k_1| + |k_2| + |l_3/k_3| + |m_4/k_4|\right) \cdot \varepsilon$$

$$|G(x,y,z)| \le \left(|l_1| + |l_2| + 1\right) \cdot \varepsilon$$

$$|H(x,y,z)| \le \left(|m_1| + |m_2| + 1\right) \cdot \varepsilon$$

$$|P(x,y,z)| \le c_0 \rho_0 \left( \sqrt{k_1^2 + l_1^2 + m_1^2} + \sqrt{k_2^2 + l_2^2 + m_2^2} \right) \cdot \varepsilon$$

Therefore $|F(x,y,z)| < \varepsilon \cdot m$ and so on for any $(x, y, z) \in R^3$. This means that $I = (F, G, H, P) \in V_0^{\varepsilon \cdot m}$. □

Let $Z_1$ be the topological function space of the set of systems $I = (F, G, H, P)$ of the form (13) obtained for a given set of constants $k_i, l_i, m_i$ $i = \overline{1,4}$ with $k_3 \cdot k_4 \ne 0$ endowed with the usual algebraic operations and the natural relative topology generated by the topology of $X_1$ [6]. A neighborhood of the origin $O$ in $Z_1$ is a set $V_0'$ of systems $I = (F, G, H, P)$ from $Z_1$ having the property that there exists $\varepsilon \in R^1$, $\varepsilon > 0$, such that if for $I = (F, G, H, P) \in Z_1$ inequalities $|F(x, y, z)| < \varepsilon$ and so on hold for any $(x, y, z) \in R^3$, then $I \in V_0'$.

The set $V_0'^\varepsilon$ defined by:

$$V_0'^\varepsilon = \{(f, G, H, P) \in Z_1 : |F(x,y,z)| < \varepsilon \text{ and so on for any } (x,y,z) \in R^3\} \quad (15)$$

is a neighborhood of the origin $O$ in $Z_1$.

As concerns the formula of representation (13) of the elements of $Z_1$ the following statement holds:

**PROPOSITION 10**. *If the constants $k_i, l_i, m_i$ $i = \overline{1,4}$ appearing in formula (13) satisfy*

$$k_1 \cdot k_2 \cdot k_3 \cdot k_4 \ne 0 \quad (16)$$

*and*

$$\Delta = \frac{c_0\rho_0}{k_3 k_4}\left[\sqrt{k_1^2 + l_1^2 + m_1^2}\,(k_2 k_3 k_4 + k_3 m_2 m_4 + k_4 l_2 l_3) + \right.$$
$$\left. \sqrt{k_2^2 + l_2^2 + m_2^2}\,(k_1 k_3 k_4 + k_3 m_1 m_4 + k_4 l_1 l_3)\right] \neq 0 \quad (17)$$

then for a system $I = (F, G, H, P) \in Z_1$ there exists a unique system $f = (f_1, f_2, f_3, f_4)$ in $Y_1$ such that equalities (13) hold. Moreover, $f_i$, $i = \overline{1,4}$ are given by:

$$f_1(\xi) = \frac{c_0\rho_0}{\Delta}\sqrt{k_2^2 + l_2^2 + m_2^2}\cdot F(\xi/k_1,0,0) + \frac{c_0\rho_0}{\Delta}\cdot\frac{l_3}{k_3}\sqrt{k_2^2 + l_2^2 + m_2^2}\cdot G(\xi/k_1,0,0) +$$
$$\frac{c_0\rho_0}{\Delta}\cdot\frac{m_4}{k_4}\sqrt{k_2^2 + l_2^2 + m_2^2}\cdot H(\xi/k_1,0,0) - \frac{1}{\Delta}\left(k_2 + \frac{m_2 m_4}{k_4} + \frac{l_2 l_3}{k_3}\right)\cdot P(\xi/k_1,0,0)$$

$$f_2(\xi) = -\frac{c_0\rho_0}{\Delta}\sqrt{k_1^2 + l_1^2 + m_1^2}\cdot F(\xi/k_1,0,0) - \frac{c_0\rho_0}{\Delta}\cdot\frac{l_3}{k_3}\sqrt{k_1^2 + l_1^2 + m_1^2}\cdot G(\xi/k_1,0,0) +$$
$$\frac{c_0\rho_0}{\Delta}\cdot\frac{m_4}{k_4}\sqrt{k_1^2 + l_1^2 + m_1^2}\cdot H(\xi/k_1,0,0) + \frac{1}{\Delta}\left(k_1 + \frac{m_1 m_4}{k_4} + \frac{l_1 l_3}{k_3}\right)\cdot P(\xi/k_1,0,0)$$

$$f_3(\xi) = \frac{c_0\rho_0}{\Delta}\left[l_1\sqrt{k_2^2 + l_2^2 + m_2^2} - l_2\sqrt{k_1^2 + l_1^2 + m_1^2}\right]\cdot F(\xi/k_1,0,0) +$$
$$\left\{\frac{c_0\rho_0}{\Delta}\left[l_1\frac{l_3}{k_3}\sqrt{k_2^2 + l_2^2 + m_2^2} - l_2\frac{l_3}{k_3}\sqrt{k_1^2 + l_1^2 + m_1^2}\right] - 1\right\}\cdot G(\xi/k_1,0,0) +$$
$$\frac{c_0\rho_0}{\Delta}\left[l_1\frac{m_4}{k_4}\sqrt{k_2^2 + l_2^2 + m_2^2} - l_2\frac{m_4}{k_4}\sqrt{k_1^2 + l_1^2 + m_1^2}\right]\cdot H(\xi/k_1,0,0) +$$
$$\frac{1}{\Delta}\left[l_2 k_1 + l_2\frac{m_1 m_4}{k_4} - l_1 k_2 - l_1\frac{m_2 m_4}{k_4}\right]\cdot P(\xi/k_1,0,0)$$

(18)

$$f_4(\xi) = \frac{c_0\rho_0}{\Delta}\left[m_1\sqrt{k_2^2 + l_2^2 + m_2^2} - m_2\sqrt{k_1^2 + l_1^2 + m_1^2}\right]\cdot F(\xi/k_1,0,0) +$$
$$\frac{c_0\rho_0}{\Delta}\left[m_1\frac{l_3}{k_3}\sqrt{k_2^2 + l_2^2 + m_2^2} - m_2\frac{l_3}{k_3}\sqrt{k_1^2 + l_1^2 + m_1^2}\right]\cdot G(\xi/k_1,0,0) +$$
$$\left\{\frac{c_0\rho_0}{\Delta}\left[m_1\frac{m_4}{k_4}\sqrt{k_2^2 + l_2^2 + m_2^2} - m_2\frac{m_4}{k_4}\sqrt{k_1^2 + l_1^2 + m_1^2}\right] - 1\right\}\cdot H(\xi/k_1,0,0) +$$
$$\frac{1}{\Delta}\left[m_2 k_1 + m_2\frac{l_1 l_3}{k_3} - m_1 k_2 - m_1\frac{l_2 l_3}{k_3}\right]\cdot P(\xi/k_1,0,0)$$

*Proof.* The uniqueness can be proven assuming the contrary. Formulas (18) can be obtained by solving (13) with respect to $f_1, f_2, f_3, f_4$. □

***Remark.*** If the constants $k_i, l_i, m_i$ $i = \overline{1,4}$ appearing in formula (13) satisfy (16) and (17) and if $I = (F, G, H, P) \in Z_1$ belongs to the neighborhood $V_0^{1\varepsilon}$, then $f = (f_1, f_2, f_3, f_4) \in Y_1$ belongs to the neighborhood $W_0^{\varepsilon \cdot M}$ where $M$ is given by:

$$M = \max \left\{ \frac{c_0 \rho_0}{|\Delta|} \sqrt{k_2^2 + l_2^2 + m_2^2} \left(1 + \left|\frac{l_3}{k_3}\right| + \left|\frac{m_4}{k_4}\right|\right) + \frac{1}{|\Delta|} \left(|k_2| + \left|\frac{l_2 l_3}{k_3}\right| + \left|\frac{m_2 m_4}{k_4}\right|\right); \right.$$

$$\frac{c_0 \rho_0}{|\Delta|} \sqrt{k_1^2 + l_1^2 + m_1^2} \left(1 + \left|\frac{l_3}{k_3}\right| + \left|\frac{m_4}{k_4}\right|\right) + \frac{1}{|\Delta|} \left(|k_1| + \left|\frac{l_1 l_3}{k_3}\right| + \left|\frac{m_1 m_4}{k_4}\right|\right);$$

$$\frac{c_0 \rho_0}{|\Delta|} \left[\sqrt{k_2^2 + l_2^2 + m_2^2} \cdot |l_1| \cdot \left(1 + \left|\frac{l_3}{k_3}\right| + \left|\frac{m_4}{k_4}\right|\right) + \sqrt{k_1^2 + l_1^2 + m_1^2} \cdot |l_2| \cdot \left(1 + \left|\frac{l_3}{k_3}\right| + \left|\frac{m_4}{k_4}\right|\right)\right]; 1; \quad (19)$$

$$\frac{1}{|\Delta|} \left(|l_2 k_1| + \left|l_2 \frac{m_1 m_4}{k_4}\right| + |l_1 k_2| + \left|l_1 \frac{m_2 m_4}{k_4}\right|\right); \frac{c_0 \rho_0}{|\Delta|} \left[\sqrt{k_2^2 + l_2^2 + m_2^2} \cdot |m_1| \cdot \left(1 + \left|\frac{l_3}{k_3}\right| + \left|\frac{m_4}{k_4}\right|\right) + \right.$$

$$\left. \sqrt{k_1^2 + l_1^2 + m_1^2} \cdot |m_2| \cdot \left(1 + \left|\frac{l_3}{k_3}\right| + \left|\frac{m_4}{k_4}\right|\right)\right]; \frac{1}{|\Delta|} \left(|m_2 k_1| + \left|m_2 \frac{l_1 l_3}{k_3}\right| + |m_1 k_2| + \left|m_1 \frac{l_2 l_3}{k_3}\right|\right) \right\}$$

The above statement can be proven by evaluating $|f_i(\xi)|$ $i = \overline{1,4}$ using equalities (18). Remark that conditions (16) and (17) are necessary conditions. For instance, if $\Delta = 0$, then it can happen that $F = G = H = P = 0$ even in the case when one or several of the functions $f_i$, $i = \overline{1,4}$, appearing in formula (13), of representation of $F, G, H, P$, are different from zero. For instance, this situation occurs when $k_1 = l_1 = m_1 = 1$; $k_2 = l_2 = m_2 = -1$; $f_1(\xi) = \sin \xi$; $f_2(\xi) = -\sin \xi$; $f_3(\xi) = f_4(\xi) = 0$ for any $\xi \in R^1$.

**DEFINITION 11.** A topological function space $Z_1$ for which relations (16) and (17) hold is the first phase space for the perturbation propagation problems (4), (8) and (6), respectively.

### 3.1. Well posedness of the instantaneous perturbation propagation problem in a phase space $Z_1$ and stability of the null solution

**PROPOSITION 12.** *In a phase space $Z_1$ the instantaneous perturbation propagation problem (4), (8) is well posed on any interval of time $[0, t_0]$ and the null solution is stable.*

*Proof.* In a phase space $Z_1$ the system of functions given by (7) is a solution of the initial value problem. The uniqueness of the solution can be obtained using similar arguments as in the *Theorem 2.16* [10] (page 89). More precisely, instead of the Sobolev space, considered in [10], for the uniqueness the Hilbert space $(L^2(R^3))^4$ has to be used as phase space and the initial perturbation has to be modified as follows: $\varphi_r(\xi) \cdot f_i(\xi)$ $i = \overline{1,4}$ with $\varphi_r(\xi) = 1$ for $|\xi| \leq r$; $\varphi_r(\xi) = 0$ for $|\xi| \geq r + 1$; $0 < \varphi_r(\xi) < 0$ for $r < |\xi| < r + 1$; $\varphi_r(\xi)$ infinitely differentiable and $r$ sufficiently large.

The stability of the null solution and implicitly the continuous dependence on the initial data on any finite interval of time $[0, t_0]$ can be proven as follows.

Consider the solution $(v_x', v_y', v_z', p')$ of the initial value problem (4), (8) given by (7) and a small real number $\varepsilon > 0$. If the instantaneous perturbation $I = (F, G, H, P)$ satisfies $I \in V_0^{\varepsilon/mM}$, then $f = (f_1, f_2, f_3, f_4)$ defining the initial perturbation $I$ satisfies $f \in W_0^{\varepsilon/m}$ and the solution of the initial value problem satisfies $(v_x', v_y', v_z', p') \in V_0^{\varepsilon}$ for any $t \geq 0$. This shows that in the phase space $Z_1$ the null solution is stable in sense of Lyapunov.□

The next statement concerns the set of the exponential growth rates of the solutions of the initial value problem (4), (8) in a phase space $Z_1$. In this paper the exponential growth rate of the solution $(v_x', v_y', v_z', p')$ is defined as:

$$v = \max \left\{ \sup_{(x,y,z)} \varlimsup_{t \to \infty} \frac{\ln |v_x'(x,y,z)|}{t}, \sup_{(x,y,z)} \varlimsup_{t \to \infty} \frac{\ln |v_y'(x,y,z)|}{t}, \sup_{(x,y,z)} \varlimsup_{t \to \infty} \frac{\ln |v_z'(x,y,z)|}{t}, \sup_{(x,y,z)} \varlimsup_{t \to \infty} \frac{\ln |p'(x,y,z)|}{t} \right\} \quad (20)$$

**PROPOSITION 13**. *The set of the exponential growth rates of the solutions of the equation (4) in a phase space $Z_1$ is the set of the real numbers. Moreover, there exist initial values for which the growth rate of the solution is equal to $+\infty$.*

*Proof.* Let $\mu$ be an arbitrary real number and $f_1(\xi) = f_2(\xi) = f_4(\xi) = 0$, $f_3(\xi) = \exp(-\mu \xi / k_3 U_0)$ for any $\xi \in R^1$. The exponential growth rate of the solution, corresponding to the initial value defined by the above system of functions and $k_3 = l_3 = m_3 = 1$, is equal to $\mu$. If $f_1(\xi) = f_2(\xi) = f_4(\xi) = 0$ and $f_3(\xi) = e^{\xi^2}$ for any $\xi \in R^1$, then the exponential growth rate of the corresponding solution is equal to $+\infty$.□

The above statement shows that in a phase space $Z_1$ the stability of the null solution coexists with solutions having strictly positive exponential growth rate. This phenomenon is impossible in a phase space having a locally convex topological vector space structure (Hilbert, Banach, normed). The explanation of the coexistence consists in the fact that in a phase space $Z_1$ the origin possess neighborhoods which are not absorbing.

Due to the existence of solutions whose exponential growth rate is equal to $+\infty$, Laplace transform with respect to $t$ does not exist for every solution. Therefore, dispersion relations cannot be derived. So, it is impossible to analyze the imaginary parts of the zero's of the dispersion relations as requires the Briggs-Bers stability criterion. It follows that in a phase space $Z_1$ the very popular Briggs-Bers stability analysis [5], [2] cannot be applied.

As concerns the exponential growth rate of a solution of (4), which corresponds to an initial data $I = (F, G, H, P) \in Z_1$ whose support (i.e. carrier) [7] is compact, the following statement holds.

**PROPOSITION 14**. *If the supports of the functions $F, G, H, P$ defining the initial data $I = (F, G, H, P)$ are compact, then the exponential growth rate of the corresponding solution is equal to zero.*

*Proof.* If the supports of the functions $F, G, H, P$ defining the initial data $I = (F, G, H, P)$ are compact, then by using (18) we deduce that there exists a finite interval $[a, b] \subset R^1$ such that the supports of the functions $f_i$, $i = \overline{1,4}$, given by (18), are included in the interval $[a,b]$. By using formula (7) and (8) we obtain that the corresponding solution $(v_x', v_y', v_z', p')$ satisfies:

$$|v_x'(x,y,z,t)| \leq m \cdot \max_{i=\overline{1,4}} \left\{ \sup_{\xi \in [a,b]} |f_i(\xi)| \right\}$$ and so on for any $(x,y,z) \in R^3$ and $t \geq 0$. Here $m$ is the constant defined by (14). Therefore, the exponential growth rate of the solution is equal to zero.□

Remark that it can happen that in a phase space $Z_1$ there exists instantaneous perturbation propagation for which velocity potential doesn't exist. For instance, if $k_i = l_i = m_i = 1$ for $i = \overline{1,4}$, $f_1(\xi) = \sin \xi$, $f_4(\xi) = \cos \xi$, $f_2(\xi) = f_3(\xi) = 0$, then for the propagation given by (7) velocity potential doesn't exist. This means that by time harmonic analysis this propagation from the phase space $Z_1$ cannot be captured.

### 3.2. Well posedness of the permanent source produced time harmonic perturbation propagation problem in a phase space $Z_1$ and stability of the null solution

**PROPOSITION 15.** *In a phase space $Z_1$ the permanent source produced time harmonic perturbation propagation problem (6) is well posed on any interval of time $[0,t_0]$ and the null solution is unstable.*

*Proof.* The system of functions $(v_x', v_y', v_z', p')$ defined by (9) is a solution of the permanent source produced time harmonic perturbation propagation problem (6) in the phase space $Z_1$ and corresponds to the amplitude $A = (F, G, H, P) \in Z_1$. The uniqueness of this solution can be obtained using similar arguments to those used in *Proposition 12*. The continuous dependence of the solution on the amplitude of perturbation on an arbitrary interval of time $[0,t_0]$ can be derived in the following way. Consider the solution $(v_x', v_y', v_z', p')$ given by (9), $\varepsilon \in R^1$, $\varepsilon > 0$ a small real number and $t_0 > 0$ an arbitrary real number ($t_0$ fixed). If the amplitude $A = (F, G, H, P) \in Z_1$ belongs to the neighborhood $V_0^{\varepsilon/mMt_0}$ of the origin $O$ in $Z_1$, then the system of functions $f = (f_1, f_2, f_3, f_4) \in Y_1$ defining the amplitude $A$ (i.e. given by (18)) belongs to the neighborhood $W_0^{\varepsilon/m}$ of the origin in $Y_1$ and the solution $(v_x', v_y', v_z', p')$ belongs to the neighborhood $V_0^{\varepsilon}$ of the origin in $Z_1$ (here $m$ is given by (14) and $M$ by (19)). This means continuous dependence of the solution on the amplitude of perturbation on the interval of time $[0,t_0]$.

In order to see the instability of the null solution in the phase space $Z_1$ let consider the amplitude $A = (F, G, H, P)$ defined by the system of constants $k_i = l_i = m_i = 1$ for $i = \overline{1,4}$ and the system of functions $f_1(\xi) = f_2(\xi) = f_4(\xi) = 0$ and $f_3(\xi) = \delta \cdot \sin \xi$ for $\xi \in R^1$ and $\delta \in R^1$, $\delta > 0$, i.e. $H = 0$; $P = 0$. According to (9), for this amplitude the solution of the time harmonic perturbation propagation problem is:

$$v_x'(x,y,z,t) = \delta \cdot h(t) \cdot \int_0^t \sin(x + y + z - U_0 t + U_0 \tau) \cdot \sin \omega_f \tau \, d\tau$$

$$v_y'(x,y,z,t) = -\delta \cdot h(t) \cdot \int_0^t \sin(x+y+z-U_0 t+U_0\tau)\cdot \sin\omega_f \tau\, d\tau$$

$v_z'(x,y,z,t) = 0$; $p'(x,y,z,t) = 0$.

For $\omega_f = U_0$ the above solution becomes:

$$v_x'(x,y,z,t) = \delta \cdot h(t) \cdot \left\{\frac{t}{2}\cdot\cos(x+y+z-U_0 t) + \frac{1}{2}\cdot\sin 2U_0 t[\sin(x+y+z-U_0 t) - \frac{1}{2}\cdot\cos(x+y+z-U_0 t)]\right\}$$

$$v_y'(x,y,z,t) = -\delta \cdot h(t) \cdot \left\{\frac{t}{2}\cdot\cos(x+y+z-U_0 t) + \frac{1}{2}\cdot\sin 2U_0 t[\sin(x+y+z-U_0 t) - \frac{1}{2}\cdot\cos(x+y+z-U_0 t)]\right\}$$

$v_z'(x,y,z,t) = 0$; $p'(x,y,z,t) = 0$.

Hence, for a given $\delta > 0$ and $x,y,z \in R^3$ there exists a sequence of moments $(t_n)$ tending to $+\infty$ such that the sequences $|v_x'(x,y,z,t_n)|$ and $|v_y'(x,y,z,t_n)|$ tend to $+\infty$. On the other hand, for an arbitrary $\varepsilon \in R^1$, $\varepsilon > 0$ if $\delta < \varepsilon$ the considered amplitude belongs to the neighborhood $V_0^\varepsilon$ of the origin in the phase space $Z_1$. In other words, the considered amplitude $A$ is so close to zero as we wish provided that $\delta$ is sufficiently small. This means that the null solution is unstable in sense of Lyapunov.□

This doesn't mean that the null solution is locally absolutely unstable in the sense defined in [8]. That is because the time harmonic perturbation considered in the proof is not localized in the space, i.e. the carrier, or support [7] of the amplitude $A$ is not compact.

When the time harmonic perturbation amplitude support is compact (i.e. the perturbation is localized in the space), then the following statement holds:

**PROPOSITION 16**. *If the supports of the functions $F,G,H,P$, defining the amplitude $A=(F,G,H,P)\in Z_1$ of the time harmonic perturbation are compact, then there exists $\delta > 0$ such that $f = (f_1, f_2, f_3, f_4)$ given by (18) satisfies $f \in W_0^\delta$ and the solution $(v_x', v_y', v_z', p')$ of the propagation problem (6) satisfies $(v_x', v_y', v_z', p') \in V_0^{M'\delta}$ where $M'$ is given by:*

$$M' = (b-a) \cdot \max \left\{ \left| \frac{k_1}{k_1 U_0 - c_0 \sqrt{k_1^2 + l_1^2 + m_1^2}} \right| + \left| \frac{k_2}{k_2 U_0 + c_0 \sqrt{k_2^2 + l_2^2 + m_2^2}} \right| + \left| \frac{l_3}{k_3^2 U_0} \right| + \left| \frac{m_4}{k_4^2 U_0} \right|; \right.$$

$$\left| \frac{l_1}{k_1 U_0 - c_0 \sqrt{k_1^2 + l_1^2 + m_1^2}} \right| + \left| \frac{l_2}{k_2 U_0 + c_0 \sqrt{k_2^2 + l_2^2 + m_2^2}} \right| + \left| \frac{1}{k_3 U_0} \right|;$$

$$\left| \frac{m_1}{k_1 U_0 - c_0 \sqrt{k_1^2 + l_1^2 + m_1^2}} \right| + \left| \frac{m_2}{k_2 U_0 + c_0 \sqrt{k_2^2 + l_2^2 + m_2^2}} \right| + \left| \frac{1}{k_4 U_0} \right|; \quad (21)$$

$$\left. \frac{c_0 \rho_0 \sqrt{k_1^2 + l_1^2 + m_1^2}}{\left| k_1 U_0 - c_0 \sqrt{k_1^2 + l_1^2 + m_1^2} \right|} + \frac{c_0 \rho_0 \sqrt{k_2^2 + l_2^2 + m_2^2}}{\left| k_2 U_0 + c_0 \sqrt{k_2^2 + l_2^2 + m_2^2} \right|} \right\}$$

and the bounded interval $[a,b]$ contains the supports of the functions $f_i$, $i = \overline{1,4}$.

*Proof.* If the supports of the functions $F, G, H, P$, defining the amplitude $A = (F, G, H, P)$ of the time harmonic perturbation are compact, then the supports of the functions $f_1, f_2, f_3, f_4$ given by (18) are also compact. It follows that there exists a finite interval $[a,b]$ which contains the supports of the functions $f_i$, $i = \overline{1,4}$. Let $\delta$ be the number defined as $\delta = \max_{i=\overline{1,4}} \left\{ \sup_{\xi \in [a,b]} |f_i(\xi)| \right\}$. By using formula (9) the functions $v_x', v_y', v_z', p'$ can be evaluated as follows:

$$|v_x'(x,y,z,t)| \leq \delta(b-a) \cdot \left| \frac{k_1}{k_1 U_0 - c_0 \sqrt{k_1^2 + l_1^2 + m_1^2}} \right| + \left| \frac{k_2}{k_2 U_0 + c_0 \sqrt{k_2^2 + l_2^2 + m_2^2}} \right| + \frac{|l_3|}{k_3^2 U_0} + \frac{|m_4|}{k_4^2 U_0}$$

$$|v_y'(x,y,z,t)| \leq \delta(b-a) \cdot \left| \frac{l_1}{k_1 U_0 - c_0 \sqrt{k_1^2 + l_1^2 + m_1^2}} \right| + \left| \frac{l_2}{k_2 U_0 + c_0 \sqrt{k_2^2 + l_2^2 + m_2^2}} \right| + \frac{1}{|m_3| U_0}$$

(22)

$$|v_z'(x,y,z,t)| \leq \delta(b-a) \cdot \left| \frac{m_1}{k_1 U_0 - c_0 \sqrt{k_1^2 + l_1^2 + m_1^2}} \right| + \left| \frac{m_2}{k_2 U_0 + c_0 \sqrt{k_2^2 + l_2^2 + m_2^2}} \right| + \frac{1}{|k_4| U_0}$$

$$|p'(x,y,z,t)| \leq \delta(b-a) \cdot \left[ \frac{c_0 \rho_0 \sqrt{k_1^2 + l_1^2 + m_1^2}}{\left| k_1 U_0 - c_0 \sqrt{k_1^2 + l_1^2 + m_1^2} \right|} + \frac{c_0 \rho_0 \sqrt{k_2^2 + l_2^2 + m_2^2}}{\left| k_2 U_0 + c_0 \sqrt{k_2^2 + l_2^2 + m_2^2} \right|} \right]$$

Hence we obtain that $(v_x', v_y', v_z', p') \in V_0^{M'\delta}$. □

**Remark**. In a phase space $Z_1$ the null solution is stable with respect to the set of time harmonic perturbation whose amplitude support is contained in a given compact set.

As concerns the exponential growth rate of the solution of the propagation problem (6) for amplitudes from a phase space $Z_1$, the following statement holds:

**PROPOSITION 17**. *For every real number $\mu$ there exists an amplitude $A = (F, G, H, P)$ in the phase space $Z_1$ such that the exponential growth rate of the solution of the propagation problem (6) is equal to $\mu$. If the support (carrier) of the amplitude $A = (F, G, H, P) \in Z_1$ is compact, then the exponential growth rate of the corresponding solution is equal to zero.*

*Proof.* For $f_1 = f_2 = f_3 = 0$ and $f_3(\xi) = \exp(-\mu\xi/k_3 U_0)$ for any $\xi \in R^1$, according to formulas (9), we have:

$$v_x'(x,y,z,t) = h(t) \cdot \frac{l_3}{k_3} \cdot e^{-\frac{\mu}{k_3 U_0}(k_3 x + l_3 y + m_3 z)} \cdot e^{\mu t} \int_0^t e^{-\mu\tau} \cdot \sin\omega_f \tau \, d\tau$$

$$v_y'(x,y,z,t) = h(t) \cdot e^{-\frac{\mu}{k_3 U_0}(k_3 x + l_3 y + m_3 z)} \cdot e^{\mu t} \int_0^t e^{-\mu\tau} \cdot \sin\omega_f \tau \, d\tau$$

$$v_z'(x,y,z,t) = 0$$

$$p'(x,y,z,t) = 0$$

Hence, the exponential growth rate of the solution is equal to $\mu$.

If the support of the amplitude $A = (F, G, H, P) \in Z_1$ is compact, then applying *Proposition 16* we obtain that the exponential growth rate of the solution is equal to zero. □

Retain that in a phase space $Z_1$ the Briggs-Bers stability analysis can not be applied for time harmonic perturbations. That is because the set of the exponential growth rates of the solutions is not bounded from above.

Retain also that in a phase space $Z_1$ it can happen that there exists permanent source produced time harmonic perturbation propagation for which velocity potential doesn't exist. For instance, for the time harmonic perturbation propagation appearing in *Proposition 16*, velocity potential doesn't exist. Therefore, the propagation in discussion can not be captured by time harmonic analysis.

## 4. Second Phase Space

Let $X_2$ be the normed function space of the set of systems $I = (F, G, H, P)$ (or $A = (F, G, H, P)$) of functions $F, G, H, P : R^3 \to R^1$ which are bounded and continuously differentiable endowed with the usual algebraic operations and the norm defined by:

$$\|I\|_{X_2} = \max\left\{\sup_{(x,y,z)}|F(x,y,z)|, \sup_{(x,y,z)}|G(x,y,z)|, \sup_{(x,y,z)}|H(x,y,z)|, \sup_{(x,y,z)}|P(x,y,z)|\right\}$$

The set $V_0^\varepsilon$ defined by $V_0^\varepsilon = \{I \in X_2 : \|I\|_{X_2} < \varepsilon\}$ is a neighborhood of the origin $O$ in $X_2$.

Let $Y_2$ be the normed function space of the set of systems $f = (f_1, f_2, f_3, f_4)$ of functions functions $f_i : R^1 \to R^1$, $i = \overline{1,4}$ which are bounded and continuously differentiable endowed with the usual algebraic operations and the norm defined by: $\|f\|_{Y_2} = \max_{i=\overline{1,4}}\left\{\sup_{\xi \in R^1}|f_i(\xi)|\right\}$.

The set $W_0^\varepsilon$ defined by $W_0^\varepsilon = \{f \in Y_2 : \|f\| < \varepsilon\}$ is a neighborhood of the origin in $Y_2$.

In case of the normed spaces $X_2, Y_2$ the following statement holds:

**PROPOSITION 18.** *If $f \in Y_2$ satisfies $\|f\|_{Y_2} < \varepsilon$ and $F, G, H, P$ are given by (13), then $I = (F, G, H, P) \in X_2$ and satisfies $\|I\|_{X_2} < m \cdot \varepsilon$ where $m$ is given by (14).*

*Proof.* It follows from *Proposition 9*. □

Let $Z_2$ be the normed function space of the systems $I = (F, G, H, P)$ of the form (13) obtained for a given set of constants $k_i, l_i, m_i$ $i = \overline{1,4}$ with $k_3 \cdot k_4 \neq 0$ endowed with the usual algebraic operations and the same norm as $X_2$, i.e.

$$\|I\|_{Z_2} = \|I\|_{X_2} = \max\left\{\sup_{(x,y,z)}|Fx,y,z|; \sup_{(x,y,z)}|G(x,y,z)|; \sup_{(x,y,z)}|H(x,y,z)|; \sup_{(x,y,z)}|P(x,y,z)|\right\}$$

The set $V_0^{1\varepsilon}$ defined by $V_0^{1\varepsilon} = \{I \in Z_2 : \|I\|_{Z_2} < \varepsilon\}$ is a neighborhood of the origin $O$ in $Z_2$.

**PROPOSITION 19.** *If the constants $k_i, l_i, m_i$, $i = \overline{1,4}$ satisfy (16) and (17), then for a system $I = (F, G, H, P) \in Z_2$ there exists a unique system $f = (f_1, f_2, f_3, f_4)$ in $Y_2$ such that equalities (13) hold and $f_i$, $i = \overline{1,4}$ are given by (18). Moreover, if $\|I\|_{Z_2} < \varepsilon$, then $\|f\|_{Y_2} < M \cdot \varepsilon$, where $M$ is given by (19).*

*Proof.* Similar to that of the *Proposition 10*. □

**DEFINITION 20.** A normed function space $Z_2$ for which relations (16) and (17) hold is the second phase space for the perturbation propagation problems (4),(8) and (6), respectively.

## 4.1. Well posedness of the instantaneous perturbation propagation problem in a phase space $Z_2$ and stability of the null solution

**PROPOSITION 21.** *In a phase space $Z_2$ the instantaneous perturbation propagation problem is well posed on any interval of time $[0, t_0]$ and the null solution is stable.*

*Proof.* The proof of the existence and uniqueness is similar as in the first phase space $Z_1$ (*Proposition 12*). The stability of the null solution and implicitly the continuous dependence on the initial data on any interval of time $[0, t_0]$ can be obtained as follows. Consider $\varepsilon > 0$ (small) and assume that the initial value $I = (F, G, H, P)$ satisfies $\|I\|_{Z_2} < \dfrac{\varepsilon}{m \cdot M}$, with $m$ given by (14) and $M$ given by (19). At the first step, according to *Proposition 19*, it follows that $\|f\|_{Y_2} < \dfrac{\varepsilon}{m}$. Hence, according to *Proposition 18*, we obtain that $\|(v_x', v_y', v_z', p)\|_{Z_2} < \varepsilon$. □

The next statement concerns the set of exponential growth rates of the solutions of the initial value problem (4),(8) in a phase space $Z_2$.

**PROPOSITION 22.** *In a phase space $Z_2$ the exponential growth rate of an arbitrary solution of the initial value problem (4),(8) is equal to zero.*

*Proof.* If $I = (F, G, H, P)$ is an arbitrary initial value from the phase space $Z_2$ and $(v_x', v_y', v_z', p)$ is the corresponding solution of (4), then the following inequality holds:
$\|(v_x', v_y', v_z', p)\|_{Z_2} \leq m \cdot \|I\|_{Z_2}$, where $m$ is given by (14).
It follows that the exponential growth rate of the solution is equal to zero. □

Due to the fact that the set of the exponential growth rate of the solutions is bounded from above by zero, for any solution the Laplace transform with respect to the time variable $t \geq 0$ exists for any complex number $\zeta \in C$ with $\text{Re}\,\zeta > 0$. However, the Briggs-Bers stability analysis can not be applied. That is because the Fourier transform with respect to the space variables doesn't exist for every solution in a phase space $Z_2$. For instance, if $k_i = l_i = m_i = 1$ for $i = \overline{1,4}$, $f_1(\xi) = f_2(\xi) = f_3(\xi) = 0$, $f_4(\xi) = \sin \xi$ for any $\xi \in R^1$, then the explicit solution is $v_x'(x,y,z,t) = \sin(x+y+z-U_0 t)$, $v_y'(x,y,z,t) = 0$, $v_z'(x,y,z,t) = \sin(x+y+z-U_0 t)$, $p'(x,y,z,t) = 0$ and has no Fourier transform with respect to the space variables $(x,y,z)$. Therefore, dispersion relations can not be derived as requires the Briggs-Beer stability criterion.

Remark also that it can happen, that in a phase space $Z_2$ there exists instantaneous perturbation propagation for which velocity potential doesn't exist. For instance, in case of the above propagation, velocity potential doesn't exist (*Proposition 6*). This means that by time harmonic analysis this propagation, from the phase space $Z_2$, can not be captured.

### *4.2. Well posedness of the permanent source produced time harmonic perturbation propagation problem in a phase space $Z_2$ and stability of the null solution*

The next statement concerns the well-posedness and the instability of the null solution in a phase space $Z_2$.

**PROPOSITION 23.** *In a phase space $Z_2$ the permanent source produced time harmonic perturbation propagation problem (6) is well posed on any finite interval of time $[0,t_0]$ and the null solution is unstable.*

*Proof.* The system of functions $(v_x', v_y', v_z', p')$ defined by (9) is a solution of the permanent source produced time harmonic perturbation propagation problem (6) in a phase space $Z_2$ and corresponds to the amplitude $A = (F,G,H,P) \in Z_2$. The uniqueness of this solution can be obtained using similar arguments to those used in *Proposition 12*. The continuous dependence of the solution on the amplitude of perturbation on an arbitrary interval of time $[0,t_0]$ can be derived as follows. Consider the solution $(v_x', v_y', v_z', p')$ given by (9), $\varepsilon \in R^1$, $\varepsilon > 0$ a small real number and $t_0 > 0$ an arbitrary real number. If the amplitude $A = (F,G,H,P) \in Z_2$ satisfies

$\|A\|_{Z_2} < \dfrac{\varepsilon}{m \cdot M \cdot t_0}$, then the system of functions $f = (f_1, f_2, f_3, f_4) \in Y_2$ given by (18) satisfies

$\|f\|_{Y_2} < \dfrac{\varepsilon}{m}$ and the solution satisfies $\|(v_x', v_y', v_z', p')\|_{Z_2} < \varepsilon$ (here $m$ is given by (14) and $M$ by (19)). This means continuous dependence of the solution on the amplitude of perturbation on the interval of time $[0,t_0]$.

In order to show the instability of the null solution in a phase space $Z_2$, consider the same amplitude as in the proof of *Proposition 15* and remark that this belongs to the phase space $Z_2$. For $\varepsilon > 0$, if $\delta < \varepsilon$, then $\|A\| < \varepsilon$ and for a given $\delta > 0$, $\delta < \varepsilon$ and $(x,y,z) \in R^3$ there exists a sequence of moments $(t_n)$ tending to $+\infty$ such that the sequences $|v'_x(x,y,z,t_n)|$, $|v'_y(x,y,z,t_n)|$ tend to $+\infty$. This means that the null solution is unstable in sense of Lyapunov. The amplitude of the perturbation is not localized in the space and the instability is a global absolute instability [8].□

When the permanent source produced time harmonic perturbation amplitude carrier is compact (i.e. localized, then the following statement holds:

**PROPOSITION 24.** *If the carriers of the functions $F,G,H,P$, defining the amplitude $A = (F,G,H,P) \in Z_2$ of the time harmonic perturbation are compact, then there exists $\delta > 0$ such that $f = (f_1,f_2,f_3,f_4)$ given by (18) satisfies $\|f\|_{Y_2} < \delta$ and the solution $(v'_x,v'_y,v'_z,p')$ of the propagation problem (6) satisfies $\|(v'_x,v'_y,v'_z,p')\|_{Z_2} < M' \cdot \delta$ where $M'$ is given by (21).*

*Proof.* Similar to that of the *Proposition 16*.□

***Corollary***. In a phase space $Z_2$ the null solution is stable with respect to the set of time harmonic perturbations whose amplitude carrier is contained in a given compact set.

As concerns the exponential growth rate of the solution of the propagation (6) for amplitudes from a phase space $Z_2$, the following statement holds.

**PROPOSITION 25.** *In a phase space $Z_2$ the exponential growth rate of an arbitrary solution of the permanent source produced time harmonic perturbation propagation problem (6) is equal to zero.*

*Proof.* An arbitrary solution $(v'_x,v'_y,v'_z,p')$ of the permanent source produced time harmonic perturbation propagation problem (6) in a phase space $Z_2$ is given by (9) with $f = (f_1,f_2,f_3,f_4)$ from $Y_2$ and satisfies the following inequality:

$$\|(v'_x,v'_y,v'_z,p')\|_{Z_2} \leq m \cdot t \cdot \|f\|_{Y_2}$$

for $t \geq 0$ and $m$ given by (14). Therefore, the exponential growth rate of the solution is equal to zero.□

Remark that in a phase space $Z_2$ the Briggs-Bers stability analysis can not be applied for time harmonic perturbations. That is because Fourier transform with respect to the spatial variable doesn't exist for every solution.

Remark also that in a phase space $Z_2$ it can happen that there exists time harmonic perturbation propagation for which velocity potential doesn't exist. For instance, if $k_i = l_i = m_i = 1$ for $i = \overline{1,4}$ and for the system of functions $f_1(\xi) = f_2(\xi) = f_4(\xi) = 0$, $f_3(\xi) = \sin\xi$, i.e. $F = f_3$, $G = -f_3$, $H = 0$, $P = 0$, then for the propagation appearing in the proof of the *Proposition 15*, velocity potential doesn't exist. Therefore, the propagation in discussion can not be captured by time harmonic analysis.

**5. Third Phase Space**

Frequently, in the published papers concerning sound attenuation in lined duct, the amplitude of the source produced permanent time harmonic perturbation is the delta distribution [4]. Since the delta distribution is a tempered distribution, we will take the third phase space $Z_3$ the space of the systems $I = (F,G,H,P)$ (or $A = (F,G,H,P)$) of tempered distributions

$F, G, H, P \in S'(R^3)$ [12] endowed with the usual algebraic operations and seminorms $q_B^{Z_3}$ defined by

$$q_B^{Z_3}(I) = \max\{q_B(F), q_B(G), q_B(H), q_B(P)\}$$

where $B$ are bounded sets in the space of the rapidly decreasing functions $S(R^3)$ defined by $q_B(F) = \max_{\varphi \in B} |F(\varphi)|$ and so on.

The set $V_0^{\varepsilon, q_B^{Z_3}}$ defined by

$$V_0^{\varepsilon, q_B^{Z_3}} = \{I \in Z_3 : q_B^{Z_3}(I) < \varepsilon\}$$

is a neighborhood of the origin $O$ in $Z_3$.

### 5.1. Well posedness of the instantaneous perturbation propagation problem and stability of the null solution in the phase space $Z_3$

If the initial value is a system of tempered distributions, $I = (F, G, H, P) \in Z_3$, then the solution of the perturbation propagation problem is a family (depending on parameter $t$) $(v'_x(t), v'_y(t), v'_z(t), (p'(t)))$ of tempered distributions which satisfy (4) for $t \geq 0$ and the initial condition

$$v'_x(0) = F \; ; \; v'_y(0) = G \; ; \; v'_z(0) = H \; ; \; p'(0) = P. \tag{23}$$

**PROPOSITION 26.** *If $(v'_x(t), v'_y(t), v'_z(t), p'(t))$ is a solution of the initial value problem (4), (23), then its Fourier transform with respect to the spatial variables $(x, y, z)$ is a family of systems of tempered distributions $(\hat{v}'_x(t), \hat{v}'_y(t), \hat{v}'_z(t), \hat{p}'(t))$, $t \geq 0$ which satisfies the following equations:*

$$\begin{aligned}
\frac{\partial \hat{v}'_x}{\partial t} &= -ikU_0 \cdot \hat{v}'_x - \frac{ik}{\rho_0} \cdot \hat{p}' \\
\frac{\partial \hat{v}'_y}{\partial t} &= -ikU_0 \cdot \hat{v}'_y - \frac{il}{\rho_0} \cdot \hat{p}' \\
\frac{\partial \hat{v}'_z}{\partial t} &= -ikU_0 \cdot \hat{v}'_z - \frac{im}{\rho_0} \cdot \hat{p}' \\
\frac{\partial \hat{p}'_x}{\partial t} &= -ik\rho_0 c_0^2 \cdot \hat{v}'_x - il\rho_0 c_0^2 \cdot \hat{v}'_y - im\rho_0 c_0^2 \cdot \hat{v}'_z - ikU_0 \cdot \hat{p}'
\end{aligned} \tag{24}$$

*and the initial condition:*

$$\hat{v}'_x(0) = \hat{F} \; ; \; \hat{v}'_y(0) = \hat{G} \; ; \; \hat{v}'_z(0) = \hat{H} \; ; \; \hat{p}'(0) = \hat{P} \tag{25}$$

*where $\alpha = (k, l, m)$ is the Fourier transform variable.*

*Proof.* By computation.□

**PROPOSITION 27.** *If* $(\hat{v}'_x(t), \hat{v}'_y(t), \hat{v}'_z(t), \hat{p}'(t))$ *is a family of tempered distributions which satisfies (24), (25), then* $\hat{v}'_x(t), \hat{v}'_y(t), \hat{v}'_z(t), \hat{p}'(t)$ *are given by the following formula:*

$$\hat{v}'_x(t) = \frac{e^{-ikU_0 t}}{\|\alpha\|^2}\left[l^2 + m^2 + \frac{k^2}{2}\left(e^{ic_0\|\alpha\|\cdot t} + e^{-ic_0\|\alpha\|\cdot t}\right)\right]\hat{F} + \frac{e^{-ikU_0 t}}{\|\alpha\|^2}\left[-kl + \frac{kl}{2}\left(e^{ic_0\|\alpha\|\cdot t} + e^{-ic_0\|\alpha\|\cdot t}\right)\right]\hat{G} +$$

$$\frac{e^{-ikU_0 t}}{\|\alpha\|^2}\left[-km + \frac{km}{2}\left(e^{ic_0\|\alpha\|\cdot t} + e^{-ic_0\|\alpha\|\cdot t}\right)\right]\hat{H} + \frac{k e^{-ikU_0 t}}{2c_0\rho_0\|\alpha\|}\left(-e^{ic_0\|\alpha\|\cdot t} + e^{-ic_0\|\alpha\|\cdot t}\right)\hat{P}$$

$$\hat{v}'_y(t) = \frac{e^{-ikU_0 t}}{\|\alpha\|^2}\left[-kl + \frac{kl}{2}\left(e^{ic_0\|\alpha\|\cdot t} + e^{-ic_0\|\alpha\|\cdot t}\right)\right]\hat{F} + \frac{e^{-ikU_0 t}}{\|\alpha\|^2}\left[k^2 + m^2 + \frac{l^2}{2}\left(e^{ic_0\|\alpha\|\cdot t} + e^{-ic_0\|\alpha\|\cdot t}\right)\right]\hat{G} +$$

$$\frac{e^{-ikU_0 t}}{\|\alpha\|^2}\left[-lm + \frac{lm}{2}\left(e^{ic_0\|\alpha\|\cdot t} + e^{-ic_0\|\alpha\|\cdot t}\right)\right]\hat{H} + \frac{l e^{-ikU_0 t}}{2c_0\rho_0\|\alpha\|}\left(-e^{ic_0\|\alpha\|\cdot t} + e^{-ic_0\|\alpha\|\cdot t}\right)\hat{P}$$

(26)

$$\hat{v}'_z(t) = \frac{e^{-ikU_0 t}}{\|\alpha\|^2}\left[-km + \frac{km}{2}\left(e^{ic_0\|\alpha\|\cdot t} + e^{-ic_0\|\alpha\|\cdot t}\right)\right]\hat{F} + \frac{e^{-ikU_0 t}}{\|\alpha\|^2}\left[-lm + \frac{lm}{2}\left(e^{ic_0\|\alpha\|\cdot t} + e^{-ic_0\|\alpha\|\cdot t}\right)\right]\hat{G} +$$

$$\frac{e^{-ikU_0 t}}{\|\alpha\|^2}\left[k^2 + l^2 + \frac{m^2}{2}\left(e^{ic_0\|\alpha\|\cdot t} + e^{-ic_0\|\alpha\|\cdot t}\right)\right]\hat{H} + \frac{m e^{-ikU_0 t}}{2c_0\rho_0\|\alpha\|}\left(-e^{ic_0\|\alpha\|\cdot t} + e^{-ic_0\|\alpha\|\cdot t}\right)\hat{P}$$

$$\hat{p}'(t) = \frac{\rho_0 c_0 k e^{-ikU_0 t}}{2\|\alpha\|}\left[-e^{ic_0\|\alpha\|\cdot t} + e^{-ic_0\|\alpha\|\cdot t}\right]\hat{F} + \frac{\rho_0 c_0 l e^{-ikU_0 t}}{2\|\alpha\|}\left[-e^{ic_0\|\alpha\|\cdot t} + e^{-ic_0\|\alpha\|\cdot t}\right]\hat{G} +$$

$$\frac{\rho_0 c_0 m e^{-ikU_0 t}}{2\|\alpha\|}\left[-e^{ic_0\|\alpha\|\cdot t} + e^{-ic_0\|\alpha\|\cdot t}\right]\hat{H} + \frac{e^{-ikU_0 t}}{2}\left[e^{ic_0\|\alpha\|\cdot t} + e^{-ic_0\|\alpha\|\cdot t}\right]\hat{P}$$

with $\|\alpha\| = (k^2 + l^2 + m^2)^{1/2}$.

*Proof.* Similar to the proof of Proposition 14 from [1]. □

It can be seen that for an arbitrary $I = (F, G, H, P) \in Z_3$ the system $(\hat{v}'_x(t), \hat{v}'_y(t), \hat{v}'_z(t), \hat{p}'(t))$ of distributions given by (26) is not necessarily a family of tempered distributions. This shows that the initial value problem in the phase space $Z_3$ is ill posed. This fact is different from the 1D case reported in [1].

### 5.2. Well posedness of the permanent source produced time harmonic perturbation propagation problem and stability of the null solution in the phase space $Z_3$

If the amplitude $A = (F, G, H, P)$ of the time harmonic source produced perturbation is a system of tempered distributions $A \in Z_3$, then the solution of the perturbation propagation problem is the imaginary part of a family (depending on parameter $t$) $(v'_x(t), v'_y(t), v'_z(t), p'(t))$ of tempered distributions which satisfies the equations:

$$\frac{\partial v'_x}{\partial t} + U_0 \cdot \frac{\partial v'_x}{\partial x} + \frac{1}{\rho_0} \cdot \frac{\partial p'}{\partial x} = h(t) \cdot e^{i\omega_f t} \cdot F$$

$$\frac{\partial v'_y}{\partial t} + U_0 \cdot \frac{\partial v'_y}{\partial x} + \frac{1}{\rho_0} \cdot \frac{\partial p'}{\partial y} = h(t) \cdot e^{i\omega_f t} \cdot G \qquad (27)$$

$$\frac{\partial v'_z}{\partial t} + U_0 \cdot \frac{\partial v'_z}{\partial x} + \frac{1}{\rho_0} \cdot \frac{\partial p'}{\partial z} = h(t) \cdot e^{i\omega_f t} \cdot H$$

$$\frac{\partial p'}{\partial t} + U_0 \cdot \frac{\partial p'}{\partial x} + \rho_0 \cdot c_0^2 \left( \frac{\partial v'_x}{\partial x} + \frac{\partial v'_y}{\partial y} + \frac{\partial v'_z}{\partial z} \right) = h(t) \cdot e^{i\omega_f t} \cdot P$$

and the condition
$$v'_x(t) = 0;\ v'_y(t) = 0;\ v'_z(t) = 0;\ p'(t) = 0 \text{ for } t \leq 0. \qquad (28)$$

If $(v'_x(t), v'_y(t), v'_z(t), p'(t))$ is a family of tempered distributions in $Z_3$ which satisfies (27), (28), then its Fourier transform with respect to the spatial variables is a family $(\hat{v}'_x(t), \hat{v}'_y(t), \hat{v}'_z(t), \hat{p}'(t))$ of tempered distributions which satisfies the equations

$$\frac{\partial \hat{v}'_x}{\partial t} = -ikU_0 \cdot \hat{v}'_x - \frac{ik}{\rho_0} \cdot \hat{p}' + h(t) \cdot e^{i\omega_f t} \cdot \hat{F}$$

$$\frac{\partial \hat{v}'_y}{\partial t} = -ikU_0 \cdot \hat{v}'_y - \frac{il}{\rho_0} \cdot \hat{p}' + h(t) \cdot e^{i\omega_f t} \cdot \hat{G}$$

$$\frac{\partial \hat{v}'_z}{\partial t} = -ikU_0 \cdot \hat{v}'_z - \frac{im}{\rho_0} \cdot \hat{p}' + h(t) \cdot e^{i\omega_f t} \cdot \hat{H} \qquad (29)$$

$$\frac{\partial \hat{p}'}{\partial t} = -ik\rho_0 c_0^2 \cdot \hat{v}'_x - il\rho_0 c_0^2 \cdot \hat{v}'_y - im\rho_0 c_0^2 \cdot \hat{v}'_z - ikU_0 \cdot \hat{p}' + h(t) \cdot e^{i\omega_f t} \cdot \hat{P}$$

and the condition:
$$\hat{v}'_x(t) = 0;\ \hat{v}'_y(t) = 0;\ \hat{v}'_z(t) = 0;\ \hat{p}'(t) = 0 \text{ for } t \leq 0. \qquad (30)$$

**PROPOSITION 28.** *If $(v'_x(t), v'_y(t), v'_z(t), p'(t))$ is a family of tempered distributions in $Z_3$ which satisfies (29), (30), then $(\hat{v}'_x(t), \hat{v}'_y(t), \hat{v}'_z(t), \hat{p}'(t))$ are given by the following formula*

$$\hat{v}'_x(t) = \frac{h(t)}{\|\alpha\|^2} \left[ \frac{(l^2+m^2)(e^{i\omega_f t} - e^{-ikU_0 t})}{i(\omega_f + kU_0)} + \frac{k^2(e^{i\omega_f t} - e^{-i(kU_0 - c_0\|\alpha\|)t})}{2i(\omega_f + kU_0 - c_0\|\alpha\|)} + \frac{k^2(e^{i\omega_f t} - e^{-i(kU_0 + c_0\|\alpha\|)t})}{2i(\omega_f + kU_0 + c_0\|\alpha\|)} \right] \hat{F} +$$

$$\frac{h(t)}{\|\alpha\|^2} \left[ -\frac{kl(e^{i\omega_f t} - e^{-ikU_0 t})}{i(\omega_f + kU_0)} + \frac{kl(e^{i\omega_f t} - e^{-i(kU_0 - c_0\|\alpha\|)t})}{2i(\omega_f + kU_0 - c_0\|\alpha\|)} + \frac{kl(e^{i\omega_f t} - e^{-i(kU_0 + c_0\|\alpha\|)t})}{2i(\omega_f + kU_0 + c_0\|\alpha\|)} \right] \hat{G} +$$

$$\frac{h(t)}{\|\alpha\|^2} \left[ -\frac{km(e^{i\omega_f t} - e^{-ikU_0 t})}{i(\omega_f + kU_0)} + \frac{km(e^{i\omega_f t} - e^{-i(kU_0 - c_0\|\alpha\|)t})}{2i(\omega_f + kU_0 - c_0\|\alpha\|)} + \frac{km(e^{i\omega_f t} - e^{-i(kU_0 + c_0\|\alpha\|)t})}{2i(\omega_f + kU_0 + c_0\|\alpha\|)} \right] \hat{H} +$$

$$\frac{h(t)k}{2c_0\rho_0\|\alpha\|} \left[ -\frac{(e^{i\omega_f t} - e^{-i(kU_0 - c_0\|\alpha\|)t})}{i(\omega_f + kU_0 - c_0\|\alpha\|)} + \frac{(e^{i\omega_f t} - e^{-i(kU_0 + c_0\|\alpha\|)t})}{i(\omega_f + kU_0 + c_0\|\alpha\|)} \right] \hat{P}$$

$$\hat{v}'_y(t) = \frac{h(t)}{\|\alpha\|^2}\left[\frac{-kl(e^{i\omega_f t} - e^{-ikU_0 t})}{i(\omega_f + kU_0)} + \frac{kl(e^{i\omega_f t} - e^{-i(kU_0 - c_0\|\alpha\|)t})}{2i(\omega_f + kU_0 - c_0\|\alpha\|)} + \frac{kl(e^{i\omega_f t} - e^{-i(kU_0 + c_0\|\alpha\|)t})}{2i(\omega_f + kU_0 + c_0\|\alpha\|)}\right]\hat{F} +$$

$$\frac{h(t)}{\|\alpha\|^2}\left[\frac{(k^2 + m^2)(e^{i\omega_f t} - e^{-ikU_0 t})}{i(\omega_f + kU_0)} + \frac{l^2(e^{i\omega_f t} - e^{-i(kU_0 - c_0\|\alpha\|)t})}{2i(\omega_f + kU_0 - c_0\|\alpha\|)} + \frac{l^2(e^{i\omega_f t} - e^{-i(kU_0 + c_0\|\alpha\|)t})}{2i(\omega_f + kU_0 + c_0\|\alpha\|)}\right]\hat{G} +$$

$$\frac{h(t)}{\|\alpha\|^2}\left[-\frac{lm(e^{i\omega_f t} - e^{-ikU_0 t})}{i(\omega_f + kU_0)} + \frac{lm(e^{i\omega_f t} - e^{-i(kU_0 - c_0\|\alpha\|)t})}{2i(\omega_f + kU_0 - c_0\|\alpha\|)} + \frac{lm(e^{i\omega_f t} - e^{-i(kU_0 + c_0\|\alpha\|)t})}{2i(\omega_f + kU_0 + c_0\|\alpha\|)}\right]\hat{H} +$$

$$\frac{h(t)l}{2c_0\rho_0\|\alpha\|}\left[-\frac{(e^{i\omega_f t} - e^{-i(kU_0 - c_0\|\alpha\|)t})}{i(\omega_f + kU_0 - c_0\|\alpha\|)} + \frac{(e^{i\omega_f t} - e^{-i(kU_0 + c_0\|\alpha\|)t})}{i(\omega_f + kU_0 + c_0\|\alpha\|)}\right]\hat{P}$$

(31)

$$\hat{v}'_z(t) = \frac{h(t)}{\|\alpha\|^2}\left[-\frac{km(e^{i\omega_f t} - e^{-ikU_0 t})}{i(\omega_f + kU_0)} + \frac{km(e^{i\omega_f t} - e^{-i(kU_0 - c_0\|\alpha\|)t})}{2i(\omega_f + kU_0 - c_0\|\alpha\|)} + \frac{km(e^{i\omega_f t} - e^{-i(kU_0 + c_0\|\alpha\|)t})}{2i(\omega_f + kU_0 + c_0\|\alpha\|)}\right]\hat{F} +$$

$$\frac{h(t)}{\|\alpha\|^2}\left[-\frac{lm(e^{i\omega_f t} - e^{-ikU_0 t})}{i(\omega_f + kU_0)} + \frac{lm(e^{i\omega_f t} - e^{-i(kU_0 - c_0\|\alpha\|)t})}{2i(\omega_f + kU_0 - c_0\|\alpha\|)} + \frac{lm(e^{i\omega_f t} - e^{-i(kU_0 + c_0\|\alpha\|)t})}{2i(\omega_f + kU_0 + c_0\|\alpha\|)}\right]\hat{G} +$$

$$\frac{h(t)}{\|\alpha\|^2}\left[\frac{(k^2 + l^2)(e^{i\omega_f t} - e^{-ikU_0 t})}{i(\omega_f + kU_0)} + \frac{m^2(e^{i\omega_f t} - e^{-i(kU_0 - c_0\|\alpha\|)t})}{2i(\omega_f + kU_0 - c_0\|\alpha\|)} + \frac{m^2(e^{i\omega_f t} - e^{-i(kU_0 + c_0\|\alpha\|)t})}{2i(\omega_f + kU_0 + c_0\|\alpha\|)}\right]\hat{H} +$$

$$\frac{h(t)m}{2c_0\rho_0\|\alpha\|}\left[-\frac{(e^{i\omega_f t} - e^{-i(kU_0 - c_0\|\alpha\|)t})}{i(\omega_f + kU_0 - c_0\|\alpha\|)} + \frac{(e^{i\omega_f t} - e^{-i(kU_0 + c_0\|\alpha\|)t})}{i(\omega_f + kU_0 + c_0\|\alpha\|)}\right]\hat{P}$$

$$\hat{p}'(t) = \frac{\rho_0 c_0 k \cdot h(t)}{2\|\alpha\|}\left[\frac{(e^{i\omega_f t} - e^{-i(kU_0 - c_0\|\alpha\|)t})}{i(\omega_f + kU_0 + c_0\|\alpha\|)} - \frac{(e^{i\omega_f t} - e^{-i(kU_0 - c_0\|\alpha\|)t})}{i(\omega_f + kU_0 - c_0\|\alpha\|)}\right]\hat{F} +$$

$$\frac{\rho_0 c_0 l \cdot h(t)}{2\|\alpha\|}\left[\frac{(e^{i\omega_f t} - e^{-i(kU_0 - c_0\|\alpha\|)t})}{i(\omega_f + kU_0 + c_0\|\alpha\|)} - \frac{(e^{i\omega_f t} - e^{-i(kU_0 - c_0\|\alpha\|)t})}{i(\omega_f + kU_0 - c_0\|\alpha\|)}\right]\hat{G} +$$

$$\frac{\rho_0 c_0 m \cdot h(t)}{2\|\alpha\|}\left[\frac{(e^{i\omega_f t} - e^{-i(kU_0 - c_0\|\alpha\|)t})}{i(\omega_f + kU_0 + c_0\|\alpha\|)} - \frac{(e^{i\omega_f t} - e^{-i(kU_0 - c_0\|\alpha\|)t})}{i(\omega_f + kU_0 - c_0\|\alpha\|)}\right]\hat{H} +$$

$$\frac{h(t)}{2}\left[\frac{(e^{i\omega_f t} - e^{-i(kU_0 + c_0\|\alpha\|)t})}{i(\omega_f + kU_0 + c_0\|\alpha\|)} + \frac{(e^{i\omega_f t} - e^{-i(kU_0 - c_0\|\alpha\|)t})}{i(\omega_f + kU_0 - c_0\|\alpha\|)}\right]\hat{P}$$

where $(k,l,m) = \alpha$ is the variable of Fourier transform and $\|\alpha\| = (k^2 + l^2 + m^2)^{1/2}$.

*Proof.* Similar to the proof of Proposition 15 from [1]. □

It can be seen that for arbitrary $A = (F, G, H, P) \in Z_3$ the system $(\hat{v}'_x(t), \hat{v}'_y(t), \hat{v}'_z(t), \hat{p}'(t))$ of distributions given by (31) is not necessarily a family of tempered distributions. This shows that the permanent source produced time harmonic perturbation propagation problem in the phase space $Z_3$ is ill posed. This fact is similar to that of the 1D case reported in [1].

## 6. Numerical Illustration

In this section we give a numerical illustration of the global absolute instability of the null solution with respect to the time harmonic perturbation of amplitude $A = (F,G,H,P)$ with $F(x,y,z) = G(x,y,z) = H(x,y,z) = \sin(x+y+z)$ and $P(x,y,z) = -c_0\rho_0\sqrt{3}\cdot\sin(x+y+z)$ with $c_0 = 345\,m/s$, $\rho_0 = 1.20\,kg/m^3$ and angular frequency $\omega_f = U_0 - c_0\sqrt{3} = -517.557\,rad/s$ in the first phase space $Z_1$ when the velocity of the uniform flow $U_0$ is equal to 80 m/s.
The propagation of the perturbation in this case is given by:

$$v_x'(x,y,z,t) = v_y'(x,y,z,t) = v_z'(x,y,z,t) = h(t)\cdot\frac{1}{2}\cdot\cos(x+y+z-\omega_f t)\cdot\left[t - \frac{1}{2\omega_f}\cdot\sin 2\omega_f t\right] +$$

$$h(t)\cdot\frac{1}{2\omega_f}\cdot\sin(x+y+z-\omega_f t)\cdot\sin^2\omega_f t$$

$$p'(x,y,z,t) = -c_0\rho_0\cdot\sqrt{3}\cdot v'_x(x,y,z,t).$$

In Figures. 1a, 1b, 1c the variations of the amplitude of the $v_x'$ velocity component are presented in an arbitrary point of the planes: $x+y+z=0$; $x+y+z=6$ and $x+y+z=100$, respectively, during the first 3600 seconds.

Figure 1. (a,b,c)

In Figures. 2a, 2b, 2c the variations of the amplitude of the sound pressure level are presented in the same places during the same period.

Figure 2. (a, b, c)

In Figure 3 the evolution of the $v_x'$ component of the perturbed velocity and sound pressure level $p'$, respectively, are presented in the points of the $Ox$ axis situated in the range $[-20,...,20]$ during the first 600 seconds.

Figure 3 (a, b)

Finally, for the illustration of the stability, when the support (carrier) of the amplitude is compact. the characteristic function $\varphi$ of the interval $[-1,1]$ was considered and the function $f_1(\xi) = \sin\xi$, which generates the amplitude, was multiplied by $\varphi$. For this case the computed evolution of $v'_x$ and $p'$ in the points of the $Ox$ axis situated in the range $[-20, 20]$ during the first 600 seconds are presented in Figures.4a and 4b, respectively.

Figure 4 (a,b)

## 7. Conclusions and comments

1. In the first phase space $Z_1$ the instantaneous perturbation propagation problem is well posed and the null solution is stable. However, the set of the exponential growth rates of the propagations is the whole real axis. There exist explicit solutions for which velocity potential does not exist and so they can not be captured by time harmonic analysis. The Briggs-Bers stability analysis is not appropriate in this context. If the carrier (support) of the instantaneous perturbation is compact, then the exponential growth rate of the propagation is equal to zero. The permanent source produced time harmonic perturbation propagation problem in the phase space $Z_1$ is also well posed, but the null solution is globally absolutely unstable. The set of the exponential growth rate of the propagations is the whole real axis. There exist explicit solutions for which velocity potential does not exist and so they can not be captured by time harmonic analysis. The Briggs-Bers stability analysis is not appropriate in this context. However, when the carrier (support) of the amplitudes of a subset of time harmonic perturbations from $Z_1$ is included in a given compact set, then the null solution is stable with respect to that subset of time harmonic perturbations and the exponential growth rate of the corresponding propagations is equal to zero.

2. In the second phase space $Z_2$ the instantaneous perturbation propagation problem is well posed and the null solution is stable. The exponential growth rate of any propagation is equal to zero. There exist explicit solutions for which velocity potential does not exist and they can not be captured by time harmonic analysis. The Briggs-Bers stability analysis is not appropriate in this context.

The permanent source produced time harmonic perturbation propagation problem in the phase space $Z_2$ is also well posed, but the null solution is globally absolutely unstable. The exponential growth rate of any propagation is equal to zero. There exist explicit solutions for which velocity potential does not exist and they can not be captured by time harmonic analysis. The Briggs-Bers stability analysis is not appropriate in this context. However, if the carrier (support) of the amplitudes of a subset of time harmonic perturbations from $Z_2$ is included in a given compact set, then the null solution is stable with respect to that subset of time harmonic perturbations.

3. In the third phase space $Z_3$ (space of tempered distributions) the instantaneous perturbation propagation problem and the source produced time harmonic perturbation propagation are ill posed.

4. In case of a real phenomena, governed by the linearized Euler equations the choice of the phase space influences in a considerable amount the obtained results.

**Acknowledgement**

This work was supported by a grant of the Romanian National Authority for Scientific Research, CNCS-UEFISCDI, project number PN-II-ID-PCE 2011-3-0171.

List of figure captions

**Figure 1.** The variation $v'_x(t)$ during the first 3600 seconds in the points of the planes:
 **a.** $x+y+z=0$; **b.** $x+y+z=6$ and **c.** $x+y+z=100$, respectively.
**Figure 2.** The variation $p'(t)$ during the first 3600 seconds in the points of the planes:
 **a.** $x+y+z=0$; **b.** $x+y+z=6$ and **c.** $x+y+z=100$, respectively.
**Figure 3.** The variations of $v'_x(x,0,0,t)$ **(a)** and $p'(x,0,0,t)$ **(b)**, respectively, for $x \in [-20, 20]$ and $t \in [0, 600]$ in case of instability.
**Figure 4**. The variations of $v'_x(x,0,0,t)$ **(a)** and $p'(x,0,0,t)$ **(b)**, respectively, for $x \in [-20, 20]$ and $t \in [0, 600]$ in case of stability.

Figure 1

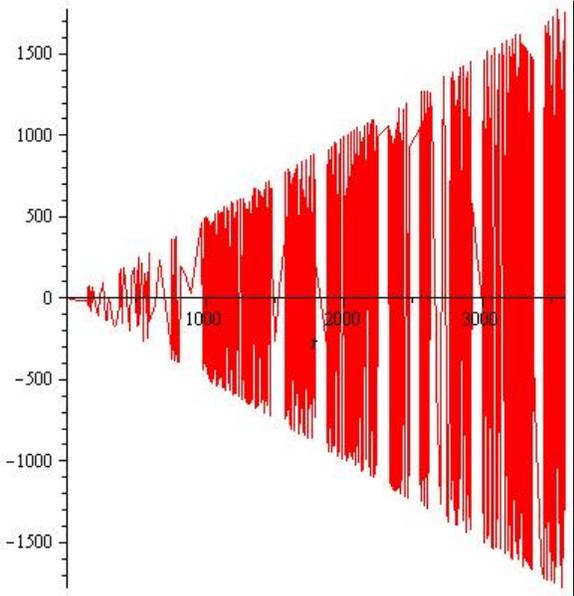

**a.**

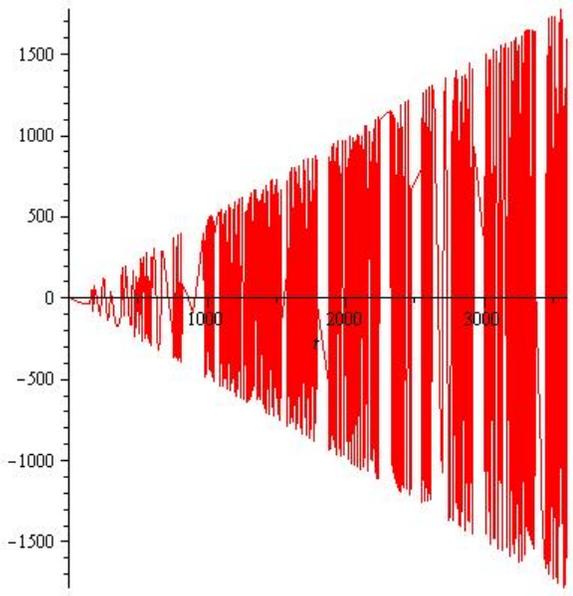

**b.**

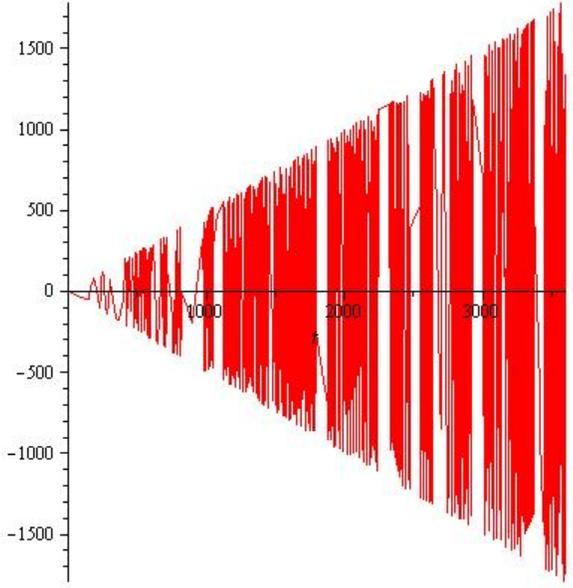

**c.**

Figure 2.

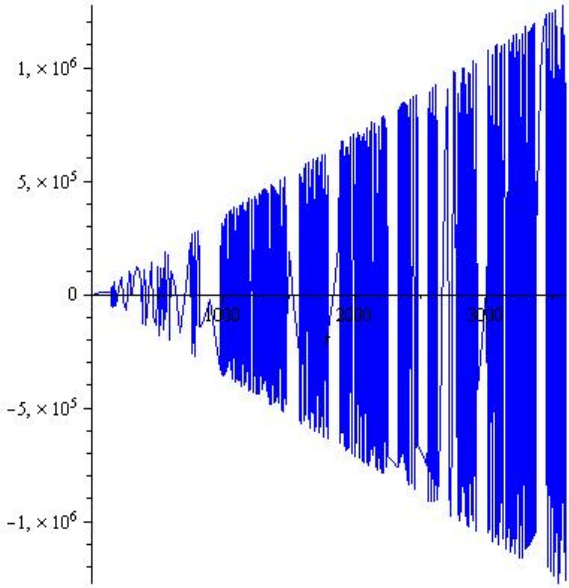

**a.**

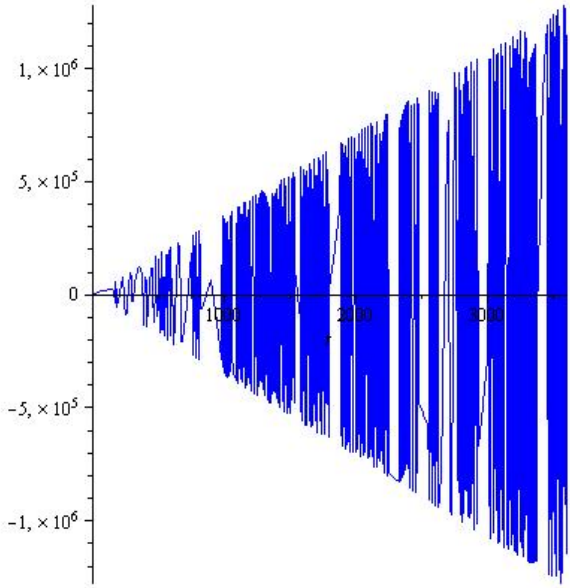

**b.**

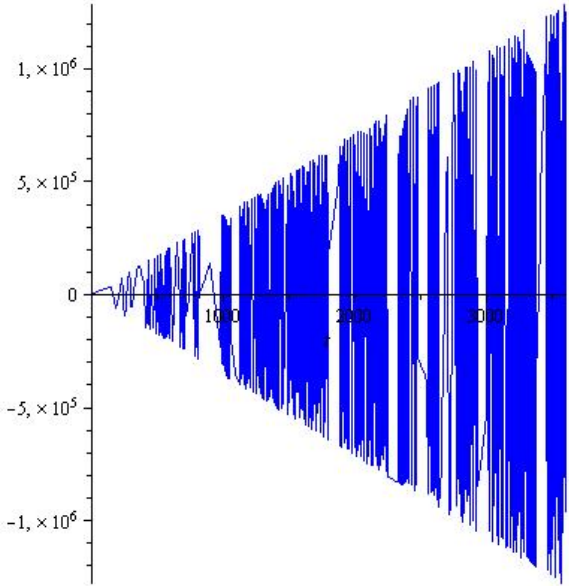

**c.**

Figure 3.

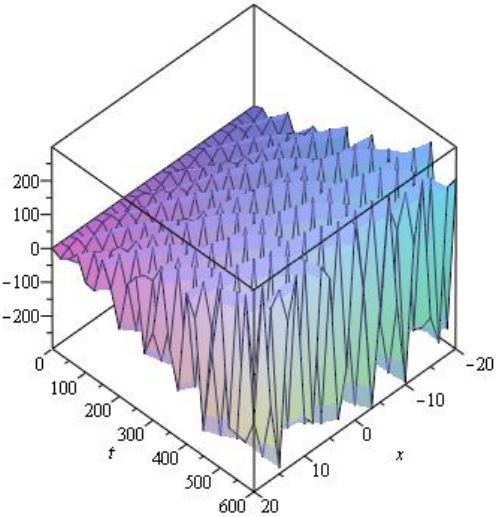 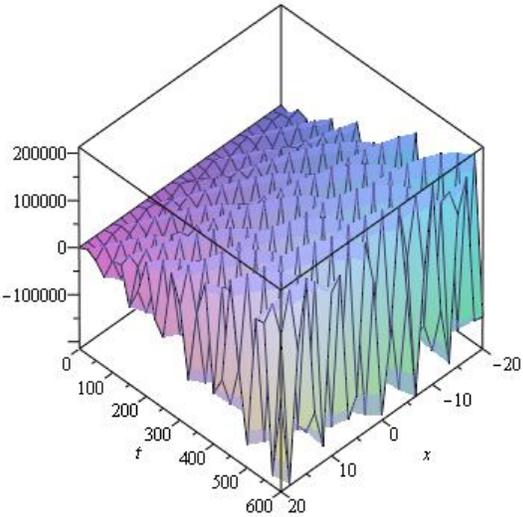

**a.** $v'_x(x,0,0,t)$          **b.** $p'(x,0,0,t)$

Figure 4

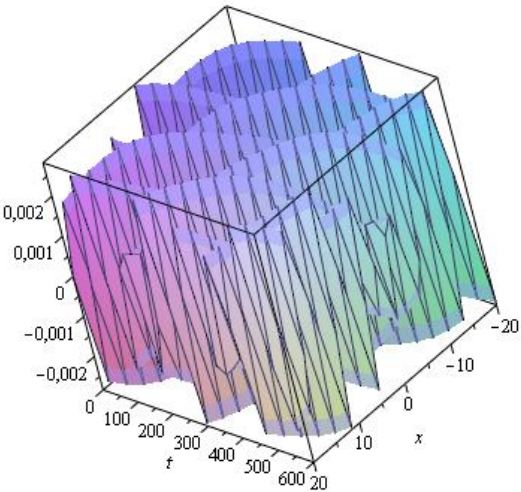 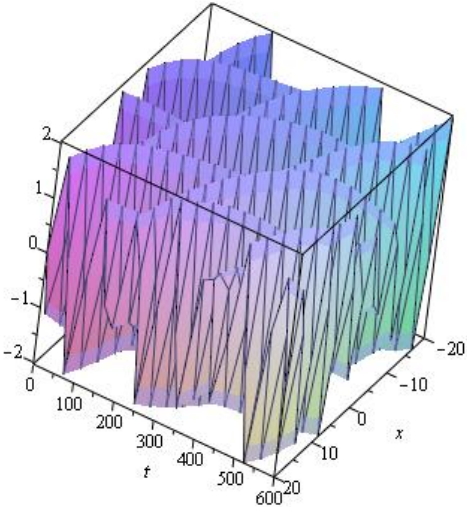

**a.** $v'_x(x,0,0,t)$          **b.** $p'(x,0,0,t)$